\documentclass[nonblindrev]{informs3} 

\usepackage{comment}
\usepackage{hyperref}
\usepackage{enumitem}

\OneAndAHalfSpacedXI 

\usepackage{natbib}
 \bibpunct[, ]{(}{)}{,}{a}{}{,}%
 \def\bibfont{\small}%
\TheoremsNumberedThrough     
\ECRepeatTheorems

\EquationsNumberedThrough    

\usepackage{amsmath}
\usepackage{float}

\theoremstyle{plain}
\newtheorem{lem}{\protect\lemmaname}
\theoremstyle{plain}

\theoremstyle{plain}

\theoremstyle{plain}

\theoremstyle{plain}
\newtheorem{algorithm}{\protect\algorithmname}
\theoremstyle{plain}
\newtheorem{thm*}{\protect\theoremname}
\theoremstyle{plain}
\newtheorem{prop*}{\protect\propositionname}

\makeatother

\usepackage[english]{babel}
\providecommand{\algorithmname}{Algorithm}
\providecommand{\corollaryname}{Corollary}
\providecommand{\definitionname}{Definition}
\providecommand{\lemmaname}{Lemma}
\providecommand{\propositionname}{Proposition}

\providecommand{\remarkname}{Remark}
\providecommand{\theoremname}{Theorem}

\begin{document}


 \RUNAUTHOR{Qu, Gao, Hinder, Ye, and Zhou}

\RUNTITLE{Optimal Diagonal Preconditioning}

\TITLE{Optimal Diagonal Preconditioning}

\ARTICLEAUTHORS{%
\AUTHOR{Zhaonan Qu}
\AFF{Department of Economics, Stanford University, \EMAIL{zhaonanq@stanford.edu}} 
\AUTHOR{Wenzhi Gao}
\AFF{School of Information Management and Engineering, Shanghai University of Finance and Economics, \EMAIL{gwz@163.shufe.edu.cn}} 
\AUTHOR{Oliver Hinder}
\AFF{Department of Industrial Engineering, University of Pittsburgh, \EMAIL{ohinder@pitt.edu}} 
\AUTHOR{Yinyu Ye}
\AFF{Department of Management Science and Engineering, Stanford University, \EMAIL{yyye@stanford.edu}}
\AUTHOR{Zhengyuan Zhou}
\AFF{Stern School of Business, New York University, \EMAIL{zzhou@stern.nyu.edu}}
} 

\ABSTRACT{%
Preconditioning has long been a staple technique in optimization, often applied to reduce the condition number of a matrix and speed up the convergence of algorithms. Although there are many popular preconditioning techniques in practice, most lack guarantees on reductions in condition number. Moreover, the degree to which we can improve over existing heuristic preconditioners remains an important practical question. In this paper, we study the problem of optimal diagonal preconditioning that achieves maximal reduction in the condition number of any full-rank matrix by scaling its rows and/or columns. We first reformulate the problem as a quasi-convex problem and provide a simple algorithm based on bisection. Then we develop an interior point algorithm with $O(\log(1/\epsilon))$ iteration complexity, where each iteration consists of a Newton update based on the Nesterov-Todd direction. Next, we specialize to one-sided optimal diagonal preconditioning problems, and demonstrate that they can be formulated as standard dual SDP problems. We then develop efficient customized solvers and study the empirical performance of our optimal diagonal preconditioning procedures through extensive experiments on large matrices. Our findings suggest that optimal diagonal preconditioners can significantly improve upon existing heuristics-based diagonal preconditioners at reducing condition numbers and speeding up iterative methods. Moreover, our implementation of customized solvers, combined with a random row/column sampling step, can find near-optimal diagonal preconditioners for matrices up to size 200,000 in reasonable time, demonstrating their practical appeal.



}%



\maketitle

%

\section{Introduction}
\label{sec:intro} 
Preconditioning, and \emph{diagonal} preconditioning in particular,
is a popular technique used in a wide range of disciplines and applications. 
Although the term ``preconditioning'' was first used in \cite{turing1948rounding}, the concept of preconditioning to improve the convergence of iterative methods traces its origins all the way back to \cite{jacobi1845ueber}, who proposed the eponymous preconditioner, which is still one of the most popular preconditioning techniques in use today. Such preconditioners are often used to speed up the convergence and improve the accuracy of iterative solvers of large linear systems, such as
as the Jacobi method \citep{saad2003iterative}, the conjugate gradient
(CG) method \citep{hestenes1952methods,barrett1994templates,concus1985block}, which is often used in Newton steps of interior point algorithms
\citep{nocedal2006numerical}, and other stationary iteration methods \citep{young1954iterative,young2014iterative} or Krylov subspace methods \citep{liesen2013krylov}.

Many first order optimization
methods also involve some form of (diagonal) preconditioning or benefit
from preconditioning. These include the ADMM method \citep{lin2018admm,takapoui2016preconditioning,boyd2011distributed}
and steepest descent methods with preconditioned gradients, such as
mirror descent with quadratic distance functions \citep{beck2003mirror},
online Newton methods, and Quasi-Newton methods \citep{broyden1970convergence,goldfarb1970family,gao2019quasi}. Intuitively, preconditioning changes the local landscape of the objective function by making its level sets more ``rounded'' so that local minima are more accessible.  

To understand the benefit of preconditioning more
precisely and to motivate the current paper, we start with the simple
setting of solving linear systems with iterative methods. Consider
the problem of solving 
\begin{equation}
\label{eq:linear-system}
Ax =b
\end{equation}
given a matrix $A\in\mathbb{R}^{m\times n}$ with full rank and $b\in\mathbb{R}^{n}$.
For large and sparse matrices $A$, iterative methods 
are very popular, in part due to parallelization. Under appropriate conditions, they often converge linearly with rates of convergence depending explicitly
on the condition number $\kappa(A)$ of $A$, defined as the ratio
of the extreme singular values of $A$:
\begin{align*}
\kappa(A) & :=\frac{\sigma_{\max}(A)}{\sigma_{\min}(A)}.
\end{align*}

\begin{table}
	\begin{centering}
	\begin{tabular}{|c|c|c|c|c|}
			\hline 
			& Jacobi & Gauss-Seidel & Steepest Descent & Conjugate Gradient\tabularnewline
			\hline 
			\hline 
			linear convergence rates &  $\frac{\kappa(A)-1}{\kappa(A)+1}$ &  $\frac{\kappa(A)-1}{\kappa(A)+1}$ &  $(\frac{\kappa(A)-1}{\kappa(A)+1})^{2}$ & $\frac{\sqrt{\kappa(A)}-1}{\sqrt{\kappa(A)}+1}$\tabularnewline
			\hline 
		\end{tabular}
		\par\end{centering}
	\caption{Rates of linear convergence of some iterative methods
		for solving the system $Ax=b$
		}
	\label{tab:rates}
\end{table}
 In Table \ref{tab:rates}, we summarize the linear convergence rates of some common iterative methods
used to solve large linear systems $Ax=b$, under their respective conditions necessary for convergence. We see that the rates are faster when the matrix $A$ has a smaller condition
number $\kappa(A)$. This phenomenon is more general.
For optimization problems, it is well known that the condition
number of the Hessian matrix of the objective plays a similar role
in the rates of convergence of first order methods. For example, for
general strongly convex objectives, gradient descent has the same
asymptotic linear convergence rate of $(\frac{\kappa(A)-1}{\kappa(A)+1})^{2}$
where $A$ is now the Hessian of the objective at optimal point $x^{\ast}$,
i.e. $A=\nabla^{2}f(x^{\ast})$. 
In order to accelerate the convergence
of iterative solvers and optimization methods, it is then essential
to ensure that the condition number of the coefficient or Hessian
matrix is small, i.e., the system is well-conditioned. Preconditioning achieves this by multiplying the
a matrix $A$ by simple matrices $D$ and $E$, such that the matrices $DA$,
$AE$, or $DAE$ have a smaller condition number than $A$. For example,
the preconditioned conjugate gradient method solves the equivalent
system
\begin{align*}
M^{-1}(Ax) & =M^{-1}b
\end{align*}
with $M$ given by the incomplete Cholesky factorization of $A$,
often resulting in $\kappa(M^{-1}A)\ll\kappa(A)$.

A large class of preconditioning methods uses matrices $D,E$ that
are \emph{diagonal}. In many applications, diagonal preconditioners are desirable, because applying them amounts to scaling the rows or columns of the matrix of interest, which does not incur much additional computational cost. In the context of solving linear systems, recovering the solution of the original system from the solution of the preconditioned system is also straightforward.

There exist many heuristics for finding
good diagonal preconditioners to achieve reductions in condition numbers.
We mention the following commonly used diagonal preconditioning procedures:
\begin{itemize}
	\item Jacobi preconditioner \citep{jacobi1845ueber}: for square matrices $A\rightarrow A(\text{diag}(A))^{-1}$ is the
	simplest diagonal preconditioner, and it works well when $A$ is diagonally
	dominant. 
	\item Matrix equilibration \citep{bradley2010algorithms, takapoui2016preconditioning}: find diagonal
	matrices $D,E$ such that columns of $DAE$ have equal $\ell_{p}$
	norms, and rows of $DAE$ have equal $\ell_{p}$ norms. A popular algorithm to achieve this is the Sinkhorn-Knopp algorithm
	\citep{sinkhorn1964relationship}. 
	\item Column (row) normalization: $A\rightarrow AD^{-1}$ where $D$ is
	the diagonal matrix of column standard deviations or column $\ell_{2}$
	norms. 
\end{itemize}

Even though the preconditioning methods discussed above are popular in practice,
they are based on heuristics and do not always guarantee reductions in the
condition number of the matrix $A$. Moreover, in many cases, an optimal diagonal preconditioner can potentially achieve a much larger reduction in condition number compared to that achieved by heuristic methods, although a detailed study of how much improvement is possible has not been conducted. This is a practically relevant issue, because it informs practitioners on whether existing popular preconditioners are sufficient for their applications. Therefore, it remains an important problem to provide an efficient and practical algorithm that provides diagonal
preconditioners which are guaranteed to achieve optimal reductions in condition numbers.


In this paper, we study the problem of
finding \emph{optimal} diagonal preconditioners for \emph{any} matrix $A$ with full rank. We are motivated by two objectives. First, we are interested in studying how much optimal diagonal preconditioners can improve over heuristic diagonal preconditioners for large systems, in terms of condition number reduction as well as speeding up iterative methods. Second, we want to develop fast and practical algorithms to find optimal diagonal preconditioners, so that when they are applied to large systems, the computational overhead to find optimal preconditioners can be potentially mitigated by overall accelerations as a result of their application. 

Our contributions are as follows. We reformulate the optimal diagonal preconditioning problem as a quasi-convex problem, and propose a bisection algorithm and  interior point algorithms with $O(\log\frac{1}{\epsilon})$
iteration complexities, that can find one-sided and two-sided diagonal preconditioners guaranteed to achieve
a condition number within $\epsilon$ of the optimal value. For one-sided preconditionings, we also provide alterantive formulations of the problem as dual SDPs, and apply customized solvers to solve large systems efficiently. We conduct extensive numerical experiments on datasets from \texttt{SuiteSparse}, \texttt{LIBSVM}, and \texttt{OPENML} to empirically evaluate the proposed procedures. Our algorithms can find near-optimal diagonal preconditioners efficiently for matrices of size up to 200,000 with the help of a row sampling scheme in the customized solver. Experiment results suggest that there is indeed a large gap between reductions in condition numbers achieved by existing popular diagonal preconditioning methods and those achieved by optimal diagonal preconditioners, suggesting to practitioners that algorithms that rely on preconditioning could be further improved. Our experiments also demonstrate that optimally-preconditioned iterative methods can outperform heuristically-preconditioned iterative methods as well as direct methods for solving linear systems.


While constructing optimal (diagonal) preconditioners incurs extra computational cost, depending on the application, it may result in reduced overall computation time. For example, if we need to solve a sequence of linear systems with identical or almost identical coefficient matrices and different right hand sides, we can use the same preconditioners repeatedly, amortizing the overhead cost of constructing them. Some applications where this is the case include: 1. solving evolution problems with implicit methods; 2. solving mildly
nonlinear problems with Newton’s method \citep{benzi2002preconditioning}; 3. solving linear systems associated with multi-scale PDE problems. Moreover, although optimal preconditioners are harder to compute for large systems, the sub-problem to which a preconditioner can be applied may not always be very large. For example, if a block-coordinate algorithm is used in an optimization procedure, each step solves a small sub-problem that may benefit from preconditioning. Likewise, the increasingly popular paradigm of Federated Learning \citep{konevcny2016federated,kairouz2021advances} requires that data be stored locally on servers, and only local parameter updates are aggregated. In this setting, each server solves a unique problem of small to moderate size, and good local preconditioners can speed up the system considerably. As another motivating example, modern interior point algorithms for solving SDPs often rely on preconditioned conjugate gradient methods in each iteration, where the system matrix remains unchanged and is usually dense \citep{andersen2011interior}. Currently, simple heuristics are used to construct diagonal preconditioners, and better preconditioners would certainly help.

Before we turn to our main problem, we briefly discuss some existing  preconditioning methods as context for the current work.
Much of the existing preconditioning literature is driven by the motivation to accelerate iterative methods, most notably the conjugate gradient, for solving large linear systems. These include the incomplete factorization methods \citep{saad2003iterative} and sparse approximate inverse methods \citep{benzi1999comparative}. The resulting preconditioners, which are generally non-diagonal, can be shown to accelerate these iterative methods. However, it is not always the case that the \emph{condition numbers} of the matrices associated with the linear system \eqref{eq:linear-system} can be provably reduced. Popular heuristics-based preconditioners such as the Jacobi preconditioner generally require diagonal dominance conditions to guarantee reductions in condition number.
In contrast, in this paper, we focus on the \emph{optimal} diagonal  preconditioning of a matrix, independent of the particular application from which it arises. This general purpose perspective allows us to develop algorithms that are guaranteed to achieve (near) optimal diagonal preconditioning for any matrix $A$ with full rank. Optimal preconditioners of this type have not been studied as much in the literature. For symmetric positive definite matrices, optimal polynomial preconditioners have been considered by \cite{cesari1937sulla,dubois1979approximating,johnson1983polynomial}. In contrast, the current work studies preconditioning for matrices of arbitrary shape, provided they have full ranks. 
The work of \cite{kaszkurewicz1995control} establishes an interesting control-theoretic interpretation of the optimal diagonal preconditioning problem for square invertible matrices, and provides a mostly theoretical algorithm by solving the associated constrained LQR problem.
The most closely related contemporary work to ours is that of \cite{jambulapati2020fast}, which considers one-sided optimal diagonal preconditioning and designs algorithms based on structured
mixed packing and covering SDPs. While \cite{jambulapati2020fast} provides theoretical guarantees on the runtime of their algorithm, it is unclear if their algorithm is practical. In contrast, we consider two-sided optimal diagonal preconditionings in addition to one-sided preconditionings, and also provide alternative SDP-based algorithms that our numerical experiments demonstrate to be effective in practice on medium to large size problems.  

Finally, we mention a few additional preconditioning methods that can serve as potential baselines with which one may compare optimal diagonal preconditioners. The first is a popular Frobenius norm minimization method \citep{benson1982iterative} to find sparse preconditioners:
\begin{align*}
    \min_{M\in \mathcal{S}}\|I-AM\|_F,
\end{align*}
where $\mathcal{S}$ is a space of sparse matrices, which we may restrict to be positive diagonal matrices, and compare the resulting diagonal preconditioners to the optimal ones, in terms of condition number reduction as well as computation cost. Another class of popular two-sided preconditioners is based on the principle of \emph{matrix equilibration} \citep{ruiz2001scaling,takapoui2016preconditioning,sinkhorn1964relationship}, among which the Ruiz scaling algorithm \citep{ruiz2001scaling} is known to perform well in practice. In our numerical experiments in Section \ref{sec:experiments}, we consider some of these baseline preconditioners when assessing the performance of our optimal preconditioners.

\section{Problem Setup and Algorithms}
\label{sec:optimal}

In this section, we turn to the two-sided optimal diagonal preconditioning
problem. For a matrix $A\in\mathbb{R}^{m\times n}$ with full rank, we want to find positive \emph{diagonal} matrices $D_1,D_2$ that solve the following minimization problem (the particular choice of exponents for $D_1,D_2$ will become clear shortly):
\begin{align}
\label{eq:two-sided-problem}
\min_{D_{1},D_{2}\succ0}\kappa(D_{1}^{1/2}AD_{2}^{-1/2})
\end{align}
There are two related \emph{one-sided} optimal diagonal preconditioning problems that seek diagonal matrices $D_1$ and $D_2$ which solve the following two distinct problems, respectively:
\begin{align}
\label{eq:left-sided-problem}
\min_{D_{1}\succ0}&\kappa(D_{1}^{1/2}A)\\
\label{eq:right-sided-problem}
\min_{D_{2}\succ0}&\kappa(AD_{2}^{-1/2})
\end{align}
 Note first that without loss of generality, we can always assume $m\geq n$, so that $A$ is a \emph{tall} matrix, since otherwise we may take the transpose of $A$ without changing the optimal solution pairs to the two-sided problem \eqref{eq:two-sided-problem}. For the one-sided problems \eqref{eq:left-sided-problem} and \eqref{eq:right-sided-problem}, taking the transpose of $A$ essentially transforms one problem to the other. Moreover, for square matrices,
i.e., $m=n$, the three optimal preconditioning problems share the
solutions, since similar square matrices
have the same eigenvalues hence the same condition numbers. Finally, note that solutions are scale-invariant, since for any optimal
pair $(D_{1},D_{2})$ of diagonal preconditioners that achieves the
minimum objective value of \eqref{eq:two-sided-problem}, $(c_1D_{1},c_2D_{2})$ for any $c_1,c_2>0$ is also
a pair of optimal diagonal preconditioners, and similarly for \eqref{eq:left-sided-problem} and \eqref{eq:right-sided-problem}. 

\textbf{Overview.} We will develop a baseline algorithm based on bisection that solves an equivalent reformulation of the two-sided problem \eqref{eq:two-sided-problem} in Section \ref{sec:bisection}, and then develop a polynomial time interior point algorithm in Section \ref{sec:potential-reduction}. This potential reduction procedure based on Newton updates with Nesterov-Todd directions can be adapted to other quasi-convex problems.  
In Section \ref{sec:DSDP}, we show that \eqref{eq:left-sided-problem} and \eqref{eq:right-sided-problem} admit reformulations as standard dual SDP problems,
to which customized solvers can be applied to obtain solutions efficiently for large matrices. We then develop specialized algorithms for the one-sided problems based on this reformulation.
In Section \ref{sec:experiments}, we study the empirical performance
of these optimal preconditioning procedures. In general, one-sided optimal diagonal preconditioners are sufficient to achieve significant reductions in condition numbers, although there are instances where two-sided preconditioners are much better. In practice, we can handle problems of size up to $m \leq 200,000$ using our customized solver with random sampling. 

\subsection{Problem Reformulation}
 We start by rewriting the two-sided optimal diagonal preconditioning problem \eqref{eq:two-sided-problem} in a more accessible form. First, note that for any matrix $A \in \mathbb{R}^{m\times n}$ with full rank and $m\geq n$, $\kappa(A)=\sqrt{\kappa(A^{T}A)}$, so we can instead minimize the condition
number of $(D_{1}^{1/2}AD_{2}^{-1/2})^{T}(D_{1}^{1/2}AD_{2}^{-1/2})=D_{2}^{-1/2}A^{T}D_{1}AD_{2}^{-1/2}$.
As the condition number is scale-invariant, we can introduce a new
positive \emph{variable} $\kappa$ (to be distinguished from the condition
number $\kappa(M)$ of a matrix) and rewrite the resulting problem
as
\begin{align*}
\min_{\kappa,D_{1},D_{2}\succ0}\kappa\\
I_n\preceq D_{2}^{-1/2}A^{T}D_{1}AD_{2}^{-1/2}\preceq\kappa I_n
\end{align*}
The objective above achieves its optimal value $\kappa^\ast:=(\min_{D_{1},D_{2}\succ0}\kappa(D_{1}^{1/2}AD_{2}^{-1/2}))^{2}$, and at any optimal solution $(D_1,D_2)$, $\lambda_{\min}(D_{2}^{-1/2}A^{T}D_{1}AD_{2}^{-1/2})=1$ and $\lambda_{\max}(D_{2}^{-1/2}A^{T}D_{1}AD_{2}^{-1/2})=\kappa^\ast$.

Next, the constraint $I_n\preceq D_{2}^{-1/2}A^{T}D_{1}AD_{2}^{-1/2}\preceq\kappa I_n$
is equivalent to $D_{2}\preceq A^{T}D_{1}A\preceq\kappa D_{2}$, and, with diagonal matrices $D_1,D_2$, we can rewrite the problem in its final reformulation as
\begin{align}
\label{eq:SDP-two-sided}
\begin{split}
    \min_{\kappa\geq0} \kappa \\
A^{T}D_{1}A\succeq D_{2}\\
\kappa D_{2}\succeq A^{T}D_{1}A \\
D_1 \succeq I_m 
\end{split}
\end{align}
where we have replaced the strict inequalities $D_1,D_2\succ0$ by equivalently requiring $\kappa \geq 0$ and $D_1 \succeq I_m$. Since $A^TD_1A\succeq A^TA\succ 0$, $D_2$ and $\kappa$ in fact need to be both strictly positive. Note that the alternative constraint $D_2 \succeq I_n$ instead of $D_1 \succeq I_m$ does not work, however, since $m\geq n$ and $D_1$ could still have zero diagonal entries even if $D_2 \succ 0$.
The problem \eqref{eq:SDP-two-sided} is also scale-invariant, and we may further impose normalization constraints to make the solution unique. This is not necessary for the purpose of the baseline algorithm we develop in this section. For the interior point algorithm we develop in Section \ref{sec:potential-reduction}, we will impose an upper bound on $D_1$ to facilitate the definition of a potential function which we will aim to reduce. 

Note also that the problem in \eqref{eq:SDP-two-sided} is a quasi-convex problem. To see this, let $d_1\in\mathbb{R}^m,d_2\in\mathbb{R}^n$ be the vector representations of the diagonal matrices $D_1,D_2$, and $d_3 \in \mathbb{R}^n$. With $d_{1i}$ denoting the $i$-th element of $d_1$, and $d_{2j},d_{3j}$ the $j$-th elements of $d_2,d_3$, \eqref{eq:SDP-two-sided} is equivalent to the following problem:
\begin{align}
\label{eq:quasi-convex}
\begin{split}
\min_{d_1,d_2,d_3\geq0} \max_j d_{3j}/d_{2j}\\
\sum_{i=1}^{m} d_{1i} a_i a_i^T - \sum_{j=1}^n d_{2j} e_j e_j^T   \succeq 0\\
\sum_{j=1}^n d_{3j} e_j e_j^T -\sum_{i=1}^{m} d_{1i} a_i a_i^T\succeq 0\\
d_1 - \mathbf{1}_m \geq 0
\end{split}
\end{align}
where vector $a_i$ is the $i$-th row of $A$ and $e_j \in \mathbb{R}^n$ is the vector whose $j$-th entry is 1 and zero otherwise. This problem is reminiscent of the Von Neumann economic growth problem \citep{neumann1945model}. 

\subsection{A Bisection Algorithm for Optimal Two-sided Diagonal Preconditioning}
\label{sec:bisection}
We start with a simple baseline algorithm for the two-sided optimal preconditioning problem \eqref{eq:SDP-two-sided} using bisection. This algorithm is easy to implement in practice, although it may be slow for larger problems. In Section \ref{sec:potential-reduction}, we will introduce an interior point algorithm with polynomial complexity, and comparison with the baseline algorithm will facilitate in understanding its improvements.

Our bisection algorithm consists of the following three steps:
\begin{enumerate}
\item Start with some large $\kappa_{0}\geq1$ such that the SDP feasibility
problem 
\begin{align*}
\min_{D_1,D_2} 0\\
A^{T}D_{1}A\succeq D_{2}\\
\kappa_{0}D_{2}\succeq A^{T}D_{1}A\\
D_1 \succeq I_m
\end{align*}
 has a solution. A choice that always works is $\kappa_{0}=\kappa(A^{T}A)$.
If we have prior information that there exists some $D_{1},D_{2}$
such that one of $\kappa(A^{T}D_{1}A)$, $\kappa(D_{2}^{-1/2}A^{T}AD_{2}^{-1/2})$,
or $\kappa(D_{2}^{-1/2}A^{T}D_{1}AD_{2}^{-1/2})$ is smaller than
$\kappa(A^{T}A)$, we can set $\kappa_{0}$ to that number as a warm
start. Popular candidates for such $D_{1},D_{2}$ could be based on
the diagonals or row/column standard deviations of $A$, or other
matrix scaling techniques such as matrix equilibriation \citep{takapoui2016preconditioning} and Ruiz scalings \citep{ruiz2001scaling}.
Set the upper bound $\overline{\kappa}=\kappa_{0}$ and lower bound
$\underline{\kappa}=1$.
\item Solve the SDP feasibility problem 
\begin{align*}
\min_{D_{1},D_{2}}0\\
A^{T}D_{1}A\succeq D_{2}\\
\frac{\overline{\kappa}+\underline{\kappa}}{2}D_{2}\succeq A^{T}D_{1}A\\
D_1 \succeq I_m
\end{align*}
If the problem is feasible, set $\frac{\overline{\kappa}+\underline{\kappa}}{2}$
to be the new \emph{upper} bound $\overline{\kappa}$. If the problem
is infeasible, set $\frac{\overline{\kappa}+\underline{\kappa}}{2}$
to be new \emph{lower} bound $\underline{\kappa}$.
\item Solve the SDP feasibility problems iteratively while updating the
upper and lower bounds, until $\overline{\kappa}-\underline{\kappa}<\epsilon$
for some tolerance parameter $\epsilon$. The final set of feasible
solutions $D_{1},D_{2}$ will provide an $\epsilon$-optimal two-sided
preconditioning.
\end{enumerate}

The algorithm can be improved via intelligent choices of the updates to lower and upper bounds by exploiting the quasi-convexity structure.
There is also a geometric interpretation of each feasibility problem
with a fixed parameter $\kappa$: it looks for $D_{1},D_{2}$ such
that the ellipsoid defined by $D_{2}$, with axes of symmetry the
coordinate axes, is inscribed in the ellipsoid defined by $A^{T}D_{1}A$,
and such that expanding $D_{2}$ by a factor of $\kappa$ will encapsulate
the ellipsoid defined by $A^{T}D_{1}A$. Such inscribing and containing
ellipsoids are used in the ellipsoid method, where the ellipsoids
are not constrained to be symmetric along coordinate axes. The geometric
interpretation also makes it clear that the optimal value $\kappa^{\ast}\geq1$
is achieved in the optimization problem. Also compare with the approach
in \cite{jambulapati2020fast}, which develops constant-factor optimal algorithms for one-sided problems based on structured
mixed packing and covering semidefinite programs. 

Note that the bisection algorithm finds an $\epsilon$-optimal solution
in $O(\log\frac{1}{\epsilon})$ iterations, where each iteration requires
solving an SDP feasibility problem. In practice, the bisection algorithm usually converges in less than 10 updates with $\epsilon = 0.01$ on most test cases of moderate size from the datasets considered in this paper (see Section \ref{sec:experiments}). 

\subsection{An Interior Point Algorithm for Optimal Diagonal Preconditioning}
\label{sec:potential-reduction}
Although the baseline algorithm introduced in the previous section is easy to implement, it suffers from sub-optimal complexity in theory and limited problem size in practice. In this section, we focus
on developing an interior point method that finds an $\epsilon$-optimal
solution in $O(\log\frac{1}{\epsilon})$ iterations, where each
iteration only consists of a Newton update step. This enjoys better
complexity than the baseline algorithm based on bisection, which requires solving an SDP feasibility problem at each iteration. Our algorithm draws inspirations from interior point
methods for linear programming and semidefinite programming
\citep{mizuno1992new,ye1992potential,ye1995neumann,nesterov1997self},
particularly \cite{nesterov1997self}
on directions of Newton updates in interior point methods for SDPs and \cite{ye1995neumann}
on the ``Von Neumann economic growth problem'' \citep{von1937uber,neumann1945model}. Our general potential reduction procedure using the Nesterov-Todd
direction is potentially applicable to other settings with quasi-convex problems, and is therefore also of independent interest.

Recall the reformulated two-sided optimal diagonal preconditioning problem in \eqref{eq:SDP-two-sided}. Due to the scale invariance nature of the problem, we impose an additional constraint on $D_1$ and develop an interior point algorithm for the following equivalent problem:
\begin{align}
\label{eq:SDP-two-sided-upper-bound}
\begin{split}
    \min_{\kappa\geq0} \kappa \\
A^{T}D_{1}A\succeq D_{2}\\
\kappa D_{2}\succeq A^{T}D_{1}A \\
\hat \kappa I_m \succeq D_1 \succeq I_m 
\end{split}
\end{align}
where $\hat{\kappa}$ is some fixed large constant. One obvious choice is $\hat{\kappa}=\kappa(A^TA)$. Following \cite{ye1995neumann}, for each fixed $\kappa\geq \kappa^\ast$, we
denote the \textbf{feasible region} of the associated SDP feasibility
problem by \textbf{
\begin{align*}
\Gamma(\kappa): & =\{D_{1},D_{2}\text{ diagonal:}\text{ }A^{T}D_{1}A\succeq D_{2};\text{ }\kappa D_{2}\succeq A^{T}D_{1}A;\text{ }\hat \kappa I_m \succeq D_{1}\succeq I_{m}\}
\end{align*}
}and define the \textbf{analytic center} of $\Gamma(\kappa)$, denoted
by $(D_{1}(\kappa),D_{2}(\kappa))$, as the \emph{unique} maximizer
of 
\begin{align*}
\max_{D_{1},D_{2}\in \Gamma(\kappa)} & \log\det(A^{T}D_{1}A-D_{2})+\log\det(\kappa D_{2}-A^{T}D_{1}A)+\log\det(D_{1}-I_{m})+\log\det(\hat \kappa I_{m}-D_1)
\end{align*}
where the objective is strictly concave because of the strict concavity
of the log determinant function. Finally, we define the \textbf{potential
function} $P(\kappa)$ as the maximum objective value above, achieved at the analytic
centers $(D_{1}(\kappa),D_{2}(\kappa))$.

As we will demonstrate, the potential function is well-defined and a suitable measure of progress towards optimality, and the iterative algorithm we design will aim to reduce the potential function by a constant amount at each step. The analytic centers will
form the central path which converges to the optimal solution as $\kappa$
is reduced together with $P(\kappa)$. See also \cite{mizuno1993adaptive,luo1998superlinear,sturm1999symmetric,burer2002solving}
for the use of potential functions and analytic centers in path-following
and potential reduction interior point algorithms for LP, SDP and
nonlinear programming problems. To facilitate presentation, we summarize our result in the following general statement and leave the detailed treatment to Appendix \ref{sec:appendix_optimal}.

\begin{theorem}
\label{thm:general-statement}
    There exists an iterative algorithm which requires a single Newton update to  $(D_1,D_2)$ at each iteration that is guaranteed to decrease the potential function $P(\kappa)$ by a constant amount. The potential is guaranteed to decrease to $O(m\log\epsilon)$ after $O(\log \frac{1}{\epsilon})$ iterations, at which point $(D_1,D_2)$ satisfies $\kappa(D_{1}^{1/2}AD_{2}^{-1/2})-\kappa^\ast = O(\epsilon)$.
\end{theorem}
Our results in this section provide us with algorithms that generate two-sided
optimal diagonal preconditioners given any fixed matrix $A$ with full rank. They can be easily adapted to the one-sided optimal preconditioning problems \eqref{eq:left-sided-problem} and \eqref{eq:right-sided-problem} as well. To our knowledge, our potential reduction algorithm is the first interior point algorithm for both the two-sided and one-sided optimal diagonal preconditioning problems. A potential concern is that the Newton updates may result in significant computation costs. In Appendix \ref{sec:first-order-methods}, we suggest some first order methods, including a first order version of the potential reduction algorithm, as alternatives to the algorithms in the main text. In practice, thanks to the sparse and low rank structure of the problems with large real datasets, our approach to find one-sided optimal diagonal preconditioners via solving dual SDPs with customized solvers, introduced in Section \ref{sec:DSDP}, is most efficient, and we may also use the procedure described in Section \ref{subsec:dual-SDP} based on these one-sided problems to obtain two-sided preconditioners. Nevertheless, the general potential reduction procedure introduced in this section can be adapted to solve other quasi-convex problems with PSD matrix variables, as is the case for problem \eqref{eq:quasi-convex}, with high accuracy when computation time is not the main constraint and bisection methods fail for large problems. 



\section{Optimal Diagonal Preconditioning via Semi-Definite Programming}

\label{sec:DSDP}

So far, we have focused on the two-side optimal diagonal preconditioning problem \begin{align*}
\min_{D_{1},D_{2}\succ0}\kappa(D_{1}^{1/2}AD_{2}^{-1/2}).
\end{align*}
In this section, we turn our attention to the one-sided optimal diagonal preconditioning problems:
\begin{align*}
\min_{D_{1}\succ0}&\kappa(D_{1}^{1/2}A)\\
\min_{D_{2}\succ0}&\kappa(AD_{2}^{-1/2}),
\end{align*} 
and propose an alternative approach to efficiently find optimal one-sided diagonal preconditioners via solving dual SDP problems. Although the
theoretical complexity of this dual SDP approach is higher than the
potential reduction algorithm we proposed in Section \ref{sec:optimal}, in practice, we are able to find optimal
diagonal preconditioners for matrices with size 5000 within reasonable time,
using an efficient customized dual SDP solver \citep{gao2022hdsdp}. Moreover, combined with a row/column sampling technique, the proposed approach can find near-optimal diagonal preconditioners for even larger (size $\geq 200,000$) problems. In Section \ref{sec:experiments}, we conduct massive numerical
experiments to explore the empirical performance of the dual SDP approach on
matrices from the \texttt{SuiteSparse} \citep{kolodziej2019suitesparse} matrix collection, as well as \texttt{LIBSVM} \citep{chang2011libsvm} and large machine learning datasets from
\texttt{OPENML} \citep{rijn2013openml}. 

\subsection{Dual SDP Formulation of Optimal Diagonal Preconditioning}
\label{subsec:dual-SDP}
Recall the following reformulation of the two-sided optimal diagonal preconditioning problem:
\begin{align*}
    \min_{\kappa\geq0} \kappa \\
A^{T}D_{1}A\succeq D_{2}\\
\kappa D_{2}\succeq A^{T}D_{1}A \\
D_1 \succeq I_m 
\end{align*}
Suppose now we are instead interested in finding optimal \emph{one-sided} diagonal preconditioners. This means setting either $D_1 = I_m$ or $D_2=I_n$ in the two-sided problem above. We will see that the resulting problems admit reformulations as standard dual SDP problems. This opens up the possibility of using efficient customized SDP solvers to find optimal diagonal preconditioners. First, setting $D_1=I_m$ gives the following optimal \emph{right} preconditioning problem:
\begin{align*}
  \min_{\kappa, D_2} \kappa \\
   A^TA \succeq D_2  \\
   \kappa D_2 \succeq A^TA  \\
   D_2 \succeq 0  
\end{align*}
where $D_2$ is restricted to be a diagonal matrix. Note that the strict lower bound $D_2 \succ 0$ is no longer necessary, since the condition that $\kappa D_2 \succeq A^TA$ will automatically enforce $\kappa > 0$ and $ D_2 \succ 0$, because we assume $A$ has full rank, and so $M:=A^TA \succ 0$.  

The key now is to do a change of variable $\tau: = 1 / \kappa$, and introduce $d:= \text{diag} (D_2 )$. Then minimization of $\kappa$ becomes \emph{maximization} of $\tau$, and with $M=A^TA$, the problem is equivalent to
\begin{align}
\label{eq:DSDP-right}
\begin{split}
      \max_{\tau, d}  \tau \\
  \text{s.t.}  \sum_{i = 1}^m e_i e_i^T d_i  \preceq M  \\
   \tau M - \sum_{i = 1}^m e_i e_i^T d_i  \preceq 0  \\
   d  \succeq 0 ,
  \end{split}
\end{align}
where $e_i$ is the $i$-th column of the identity matrix. Note that problem \eqref{eq:DSDP-right} is in standard dual SDP form (a.k.a. inequality form) \citep{luenberger1984linear}, and can be solved efficiently by customized dual SDP solvers. As we demonstrate in experiments in Section \ref{sec:experiments}, an important idea that allows us to further speed up this approach in practice for very tall matrices ($m$ very large) is to randomly \emph{sub-sample} the rows of $A$ before forming $M=A^TA$ in \eqref{eq:DSDP-right}.

Similarly, the optimal \emph{left} diagonal preconditioning problem 
\begin{align*}
    \min_{\kappa\geq0} \kappa \\
A^{T}D_{1}A\succeq I_n\\
\kappa I_n\succeq A^{T}D_{1}A \\
D_1 \succeq 0 
\end{align*}
can be transformed into a standard dual SDP with the change of variables $\tau = 1/\kappa$, $D_1/\kappa = \text{diag}(d)$:
\begin{align}
\label{eq:DSDP-left}
\begin{split}
     \max_{\tau,d}\tau\\
\text{s.t.}\sum_{i}A_{i}A_{i}^{T}d_{i} \succeq\tau I_n	\\
I_n\succeq \sum_{i}A_{i}A_{i}^{T}d_{i} \\
d\succeq0
\end{split}	
\end{align}
where we have used the identity that $A^{T}D_1A=\sum_{i}A_{i}A_{i}^{T}(D_{1})_{ii}$, where $A_{i}$ is the $i$-th row of $A$.
\begin{theorem}
The one-sided optimal diagonal preconditioning problems \eqref{eq:left-sided-problem} and \eqref{eq:right-sided-problem} are equivalent to the standard dual SDP problems \eqref{eq:DSDP-left} and \eqref{eq:DSDP-right}, respectively.
\end{theorem}

Importantly, the dual SDP formulations \eqref{eq:DSDP-right} and \eqref{eq:DSDP-left} of the one-sided problems
enjoy low-rank or sparse data matrices, and are therefore particularly amenable
to efficient customized dual SDP solvers such as \texttt{DSDP} \citep{benson2008algorithm} and \texttt{HDSDP} \citep{gao2022hdsdp}.
In Section \ref{subsec:hdsdp}, we give a brief introduction to the \texttt{HSDSP} solver, and discuss
how this dual SDP solver can solve the one-sided optimal diagonal
preconditioning problems efficiently by combining randomization and
customized numerical linear algebra.

Before we move on to discuss solving large one-sided preconditioning problems in practice, we remark on the two-sided problem. One may wonder whether the two-sided optimal diagonal preconditioning problem \eqref{eq:SDP-two-sided} also admits a reformulation as a dual SDP problem. However, because now variables $D_1$ and $D_2$ appear simultaneously on different sides of the constraint $\kappa D_{2}\succeq A^{T}D_{1}A$, it is not possible to use the same change of variable $\tau=1/\kappa$ to circumvent the product of variables as before. The equivalent quasi-convex problem \eqref{eq:quasi-convex} suggests that it may be impossible to cast the problem as an SDP. However, there exists a simple iterative procedure to obtain two-sided optimal diagonal preconditioners based on solving one-sided problems. More precisely, let $A_0 := A$. At the $k$-th step for odd $k$, we solve the left-sided problem $\min_{D_{1}\succ0}\kappa(D_{1}^{1/2}A_{k-1})$ by solving the equivalent SDP \eqref{eq:DSDP-left}, and update $A_{k}:=D_{1}^{1/2}A_{k-1}$, where $D_{1}$ is the optimal left diagonal preconditioner we just obtained. At the $k$-th step for even $k$, we solve the right-sided problem $\min_{D_{2}\succ0}\kappa(A_{k-1}D_{2}^{-1/2})$ by solving the equivalent SDP \eqref{eq:DSDP-right}, and then update $A_{k}:=A_{k-1}D_{2}^{-1/2}$.
We can then repeat the process, solving one-sided problems for the updated scaled matrix $A_k$ obtained from the previous step.
It is obvious that any fixed point of this iterative procedure is a solution to the two-sided problem. One should expect the procedure to converge, although we do not formally prove it in this paper.

\subsection{HDSDP: A Customized Dual SDP Solver}
\label{subsec:hdsdp}

In this section, we introduce the dual semi-definite programming solver
\texttt{HDSDP} \citep{gao2022hdsdp} and its internal customization mechanisms for the
optimal diagonal preconditioning problem. \texttt{HDSDP} implements a dual
potential reduction interior point method that exploits both low-rankness and
sparsity and is specifically optimized for the above preconditioning problem via
randomization.

In practice, when \texttt{HDSDP} receives $A \in \mathbb{R}^{m \times
n}$, it randomly samples $\tilde{m}$ rows without replacement and constructs a
submatrix $\tilde{A} \in \mathbb{R}^{\tilde{m} \times n}$. Then it
computes $\tilde{M} = \tilde{A}^T \tilde{A}$ and solves the right preconditioning
problems using the dual SDP \eqref{eq:DSDP-right}. Using special data
structures, the flops in each main iteration of the interior point method is
reduced to $O (n^3)$ and if $m \gg n$, we
compute the preconditioner with $O (n^3 + \tilde{m}n^2)$ complexity. We note that
when $A$ is dense and $m$ is large, computing and factorizing $A^T A$ becomes
prohibitive and the only viable method to solve $A^T A x = A^Tb$ is using
iterative solvers. By applying optimal diagonal right preconditioners, we can significantly reduce the computation time of such iterative solvers for the modified system.

Random row sampling is an important component of our solver at finding optimal preconditioners efficiently in practice, and as the empirical results in Section \ref{sec:experiments} suggest, we often only need
$\tilde{m} = O(\log{m})$ to give a high-quality preconditioner, and the overall complexity is then $O (n^3 + n^2 \log{m})$. Row sampling works particularly well when the rows of $A$ are generated from a common statistical distribution.
Note, however, that the row sampling technique only applies when finding right preconditioners.  For left preconditioning problems, we can similarly employ column sampling, although this may not yield as significant reductions in computation time since $m\geq n$.

Lastly, we note that the \texttt{HDSDP} solver we used can be modified to be a first-order type solver, where each iteration is solved by the diagonal-preconditioned conjugate gradient
method (here we just use heuristic diagonal preconditioners). Since the SDPs from large real datasets we consider are sparse, the associated SDP problems are generally handled better by 
the Cholesky method.

\section{Numerical Experiments}
\label{sec:experiments}

To empirically validate the performance of the optimal two-sided preconditioners discussed in Section \ref{sec:optimal} and the optimal one-sided preconditioners via the
dual SDP approach in Section \ref{sec:DSDP}, we conduct
extensive numerical experiments over \texttt{SuiteSparse},
\texttt{LIBSVM} and large machine learning datasets from
\texttt{OPENML}.\footnote{Code and detailed results can be found in the Github repository at \href{https://github.com/Gwzwpxz/opt\_dpcond}{https://github.com/Gwzwpxz/opt\_dpcond}} For the sake of brevity, given a matrix $A$, recall that we refer to $D_1^{1/2}A$, $AD_2^{-1/2}$ and $D_1^{1/2}AD_2^{-1/2}$ as the left, right, and two-sided preconditionings, respectively.

\subsection{Experimental Setup}

We first describe our experiment setup and implementation details.

\paragraph{Dataset.}
We use matrices from the following sources:
\begin{itemize}
  \item Tim Davis \texttt{SuiteSparse} Dataset \citep{kolodziej2019suitesparse}
  
  \item \texttt{LIBSVM} \citep{chang2011libsvm}
  
  \item \texttt{OPENML} Machine Learning Dataset Repository  \citep{rijn2013openml}

 \item Randomly generated dataset
\end{itemize}

\paragraph{Data Generation, Selection, and Processing.}
For \texttt{LIBSVM} and \texttt{OPENML}, we select their
regression datasets. Given a matrix $A$, when doing right preconditioning, if $M = A^T A$ is available, we compute $\kappa (M)$, and if it is too large, we regularize it setting $M \leftarrow M + \varepsilon I$ with $\varepsilon$ large enough such that
$\kappa (M) \leq 10^8$. This is to ensure that the associated dual SDPs are solvable by \texttt{HDSDP}. No processing is done for left or two-sided preconditionings, and if the solver fails due to ill-conditioning, we drop the matrix. For randomly generated data, we obtain their nonzero entries from $\mathcal{N} (0, 1)$.

\paragraph{Matrix Statistics.} We take matrices of different sizes from the three external
sources and record their statistics below. Details on \texttt{SuiteSparse} can be found in the Github repository.

\begin{itemize}
  \item \texttt{SuiteSparse}. We choose 391 matrices with $m, n \leq 1000$.
  
  \item \texttt{LIBSVM}. We summarize the sizes of the matrices below.
  
  \begin{table}[h]
  \centering
    \begin{tabular}{cccccc}
      \hline
      Matrix & $m$ & $n$ & Matrix & $m$ & $n$\\
            \hline
      YearPredictionMSD & 463715 & 90 & eunite2001.txt & 336 & 16\\
      YearPredictionMSD.t & 51630 & 90 & housing\_scale.txt & 506 & 13\\
      abalone\_scale.txt & 4177 & 8 & mg\_scale.txt & 1385 & 6\\
      bodyfat\_scale.txt & 252 & 14 & mpg\_scale.txt & 392 & 7\\
      cadata.txt & 20640 & 8 & pyrim\_scale.txt & 74 & 27\\
      cpusmall\_scale.txt & 8192 & 12 & space\_ga\_scale.txt & 3107 & 6\\
      eunite2001.t & 31 & 16 & triazines\_scale.txt & 186 & 60\\
      \hline
    \end{tabular}
    \caption{\texttt{LIBSVM} datasets}
  \end{table}
  \item \texttt{OPENML.} We selected the following three datasets summarized in Table \ref{tab:OPENML}.
  \begin{table}[h]
  \centering
    \begin{tabular}{ccc}
      \hline
      Matrix & $m$ & $n$\\
      \hline
      MNIST & 70000 & 784\\
      credit & 284807 & 29\\
      IMDB & 120919 & 1001\\
      \hline
    \end{tabular}
    \caption{\texttt{OPENML} datasets}
    \label{tab:OPENML}
  \end{table}
\end{itemize}

\paragraph{Subproblem Solution.}
In the experiments for right preconditioning, we adopt \texttt{HDSDP} with and without random sampling for the SDPs; for left preconditioning, we adopt \texttt{HDSDP} without random sampling;
to compute the two-sided preconditioner using bisection, we employ the better option between \texttt{Sedumi} and \texttt{SDPT3} for the feasibility problems.

\paragraph{Evaluation and Performance Metrics.} We report different statistics and metrics for the three types of preconditioning problems.

\textbf{R}ight preconditioning:

\begin{enumerate}[label={R\arabic*:},leftmargin=40pt]
  \item $\kappa (M)$: Condition number of $M=A^TA$
  
  \item $\kappa (D_2^{- 1 / 2} M D_2^{- 1 / 2})$: Condition number after optimal right preconditioning
  
  \item $\kappa (D_M^{- 1 / 2} M D_M^{- 1 / 2})$ : Condition number after Jacobi preconditioning $D_M=\text{diag}(M)$
  
  \item Improvement: Relative improvement in condition number $
  \frac{\kappa(M)}{\kappa (D_2^{- 1 / 2} M D_2^{- 1 / 2})}$
  \item Time: The CPU time spent solving the SDP problem
\end{enumerate}

\textbf{L}eft preconditioning:
\begin{enumerate}[label={L\arabic*:},leftmargin=40pt]
  \item $\kappa (M)$: Condition number of $M=A^TA$
  \item $\kappa (A^T D_1 A)$: Condition number after optimal left preconditioning  
  \item Improvement: Relative improvement in condition number $
  \frac{\kappa(M)}{\kappa (A^T D_1A)}$
  \item Time: The CPU time spent solving the SDP problem
\end{enumerate}

\textbf{T}wo-sided preconditioning:
\begin{enumerate}[label={T\arabic*:}, leftmargin=40pt]
  \item $\kappa (M)$: Condition number of $M=A^TA$
  \item $\kappa(D_2^{- 1 / 2}  A^T D_1A D_2^{- 1 / 2}) $: Condition number after optimal two-sided preconditioning  
  \item Improvement: Relative improvement in condition number $
  \frac{\kappa(M)}{\kappa(D_2^{- 1 / 2}  A^T D_1A D_2^{- 1 / 2})}$
  \item $\kappa(D_{R_1}^{- 1 / 2}  A^T D_{R_2}A D_{R_1}^{- 1 / 2}) $: Condition number after $\ell_\infty$-Ruiz preconditioning \citep{ruiz2001scaling}
\end{enumerate}

For convenience, we sometimes use Cbef, Caft to denote condition number of a matrix before and after preconditioning.

\paragraph{Conjugate Gradient Solver.}
We adopt the \texttt{MATLAB pcg} function to test the efficacy of different preconditioners at speeding up CG. The relative tolerance is set to $10^{-6}$, starting from the all-zero vector, and the right-hand-side vector is generated with a fixed random seed. When comparing the conjugate gradient method against the direct solvers, we also adopt a relative tolerance of $10^{-6}$.

\paragraph{Environment.} All experiments in this paper are carried out using \texttt{MATLAB} on
Mac Mini with Apple Silicon and 16 GB memory. The \texttt{HDSDP} solver is
implemented in \texttt{ANSI C}.

\subsection{Experiment Results for Right Preconditioning}

Our experiments for right preconditioning are into three parts. In the first part, we consider exact optimal preconditioning without randomly sampling rows of a matrix to explore the advantage of optimal preconditioners over preconditioning by the diagonal of a matrix (the Jacobi preconditioner).
In particular, we give an example where diagonal
preconditioning fails completely, but our optimal right preconditioner greatly reduces the condition number. In the
second part, we validate our randomization scheme using large realistic
regression datasets. In particular, we demonstrate that solving the dual SDP based on a sub-sampled matrix generally produces high quality approximate preconditioners for the original matrix, which has important practical implications. In the third part, we explore some cases where the time to
set up $A^T A$ is dominant, but our dual SDP approach based on random sampling can still find near-optimal preconditioners efficiently.

\paragraph{Preconditioning Evaluation.}
In this part of the experiment, we use $A$ to construct the pre-conditioner without randomly sampling its rows. The detailed results for each matrix are available in the Github repository. Here we record the summary statistics for \texttt{SuiteSparse} in Table \ref{tab:suite-right}. Note that ``improvement'' is defined as $
  \frac{\kappa(M)}{\kappa (D_2^{- 1 / 2} M D_2^{- 1 / 2})}$, so that, e.g., an improvement equal to 2 corresponds to a 2-fold improvement in the condition number.
\begin{table}[h]
\centering
  \begin{tabular}{c|c}
    \hline
    Improvement & Number\\
    \hline
    $\geq$5.00 & 121\\
    $\geq$2.00 & 190\\
    $\geq$1.25 & 261\\
    \hline
  \end{tabular}\quad\begin{tabular}{r|c}
    \hline
    Median imp. & 2.17 \\
    \hline
    Median imp. over Jacobi & 1.49 \\
    \hline
    Average time & 1.29\\
    \hline
  \end{tabular}
  \caption{Summary for right preconditioning on the \texttt{SuiteSparse} dataset}
  \label{tab:suite-right}
\end{table}

\begin{figure}
    \centering
    \includegraphics[scale=0.175]{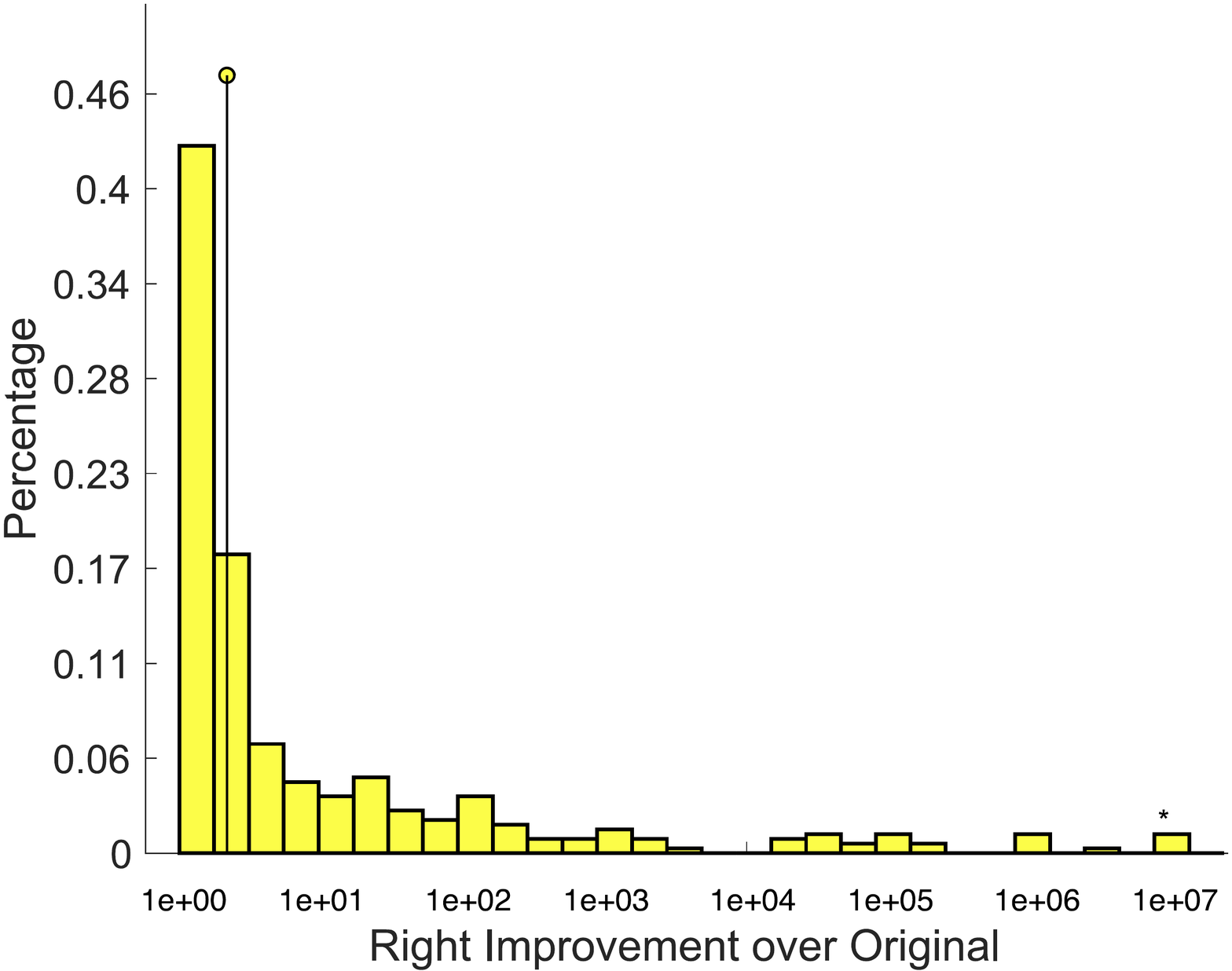}
    \includegraphics[scale=0.175]{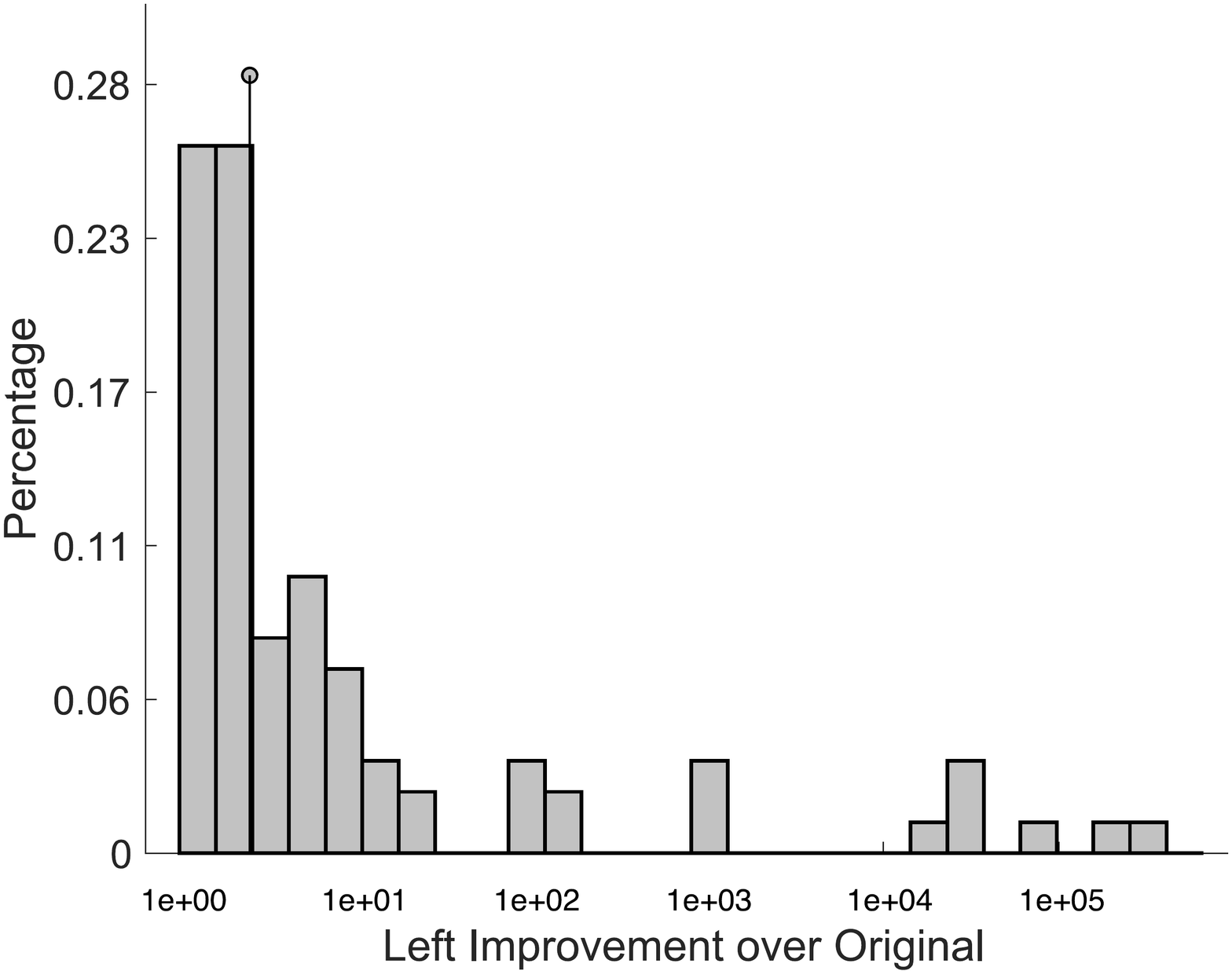}
    \includegraphics[scale=0.175]{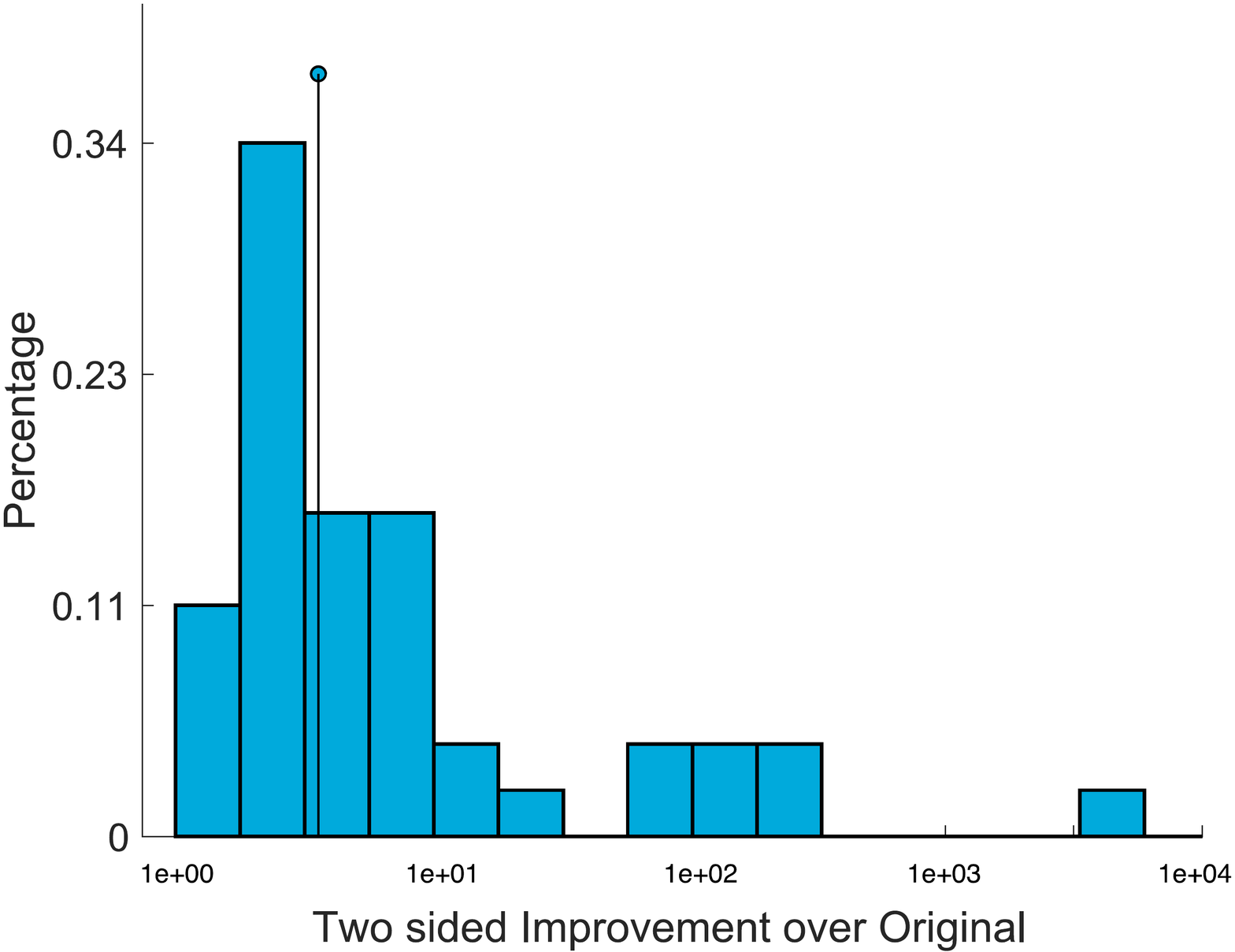}
    \caption{Distribution of logarithm of condition number improvements achieved by optimal right, left, and two-sided preconditioners on \texttt{SuiteSparse} matrices, respectively}
    \label{fig:reduction-pdf}
\end{figure}

\begin{figure}
    \centering
    \includegraphics[scale=0.18]{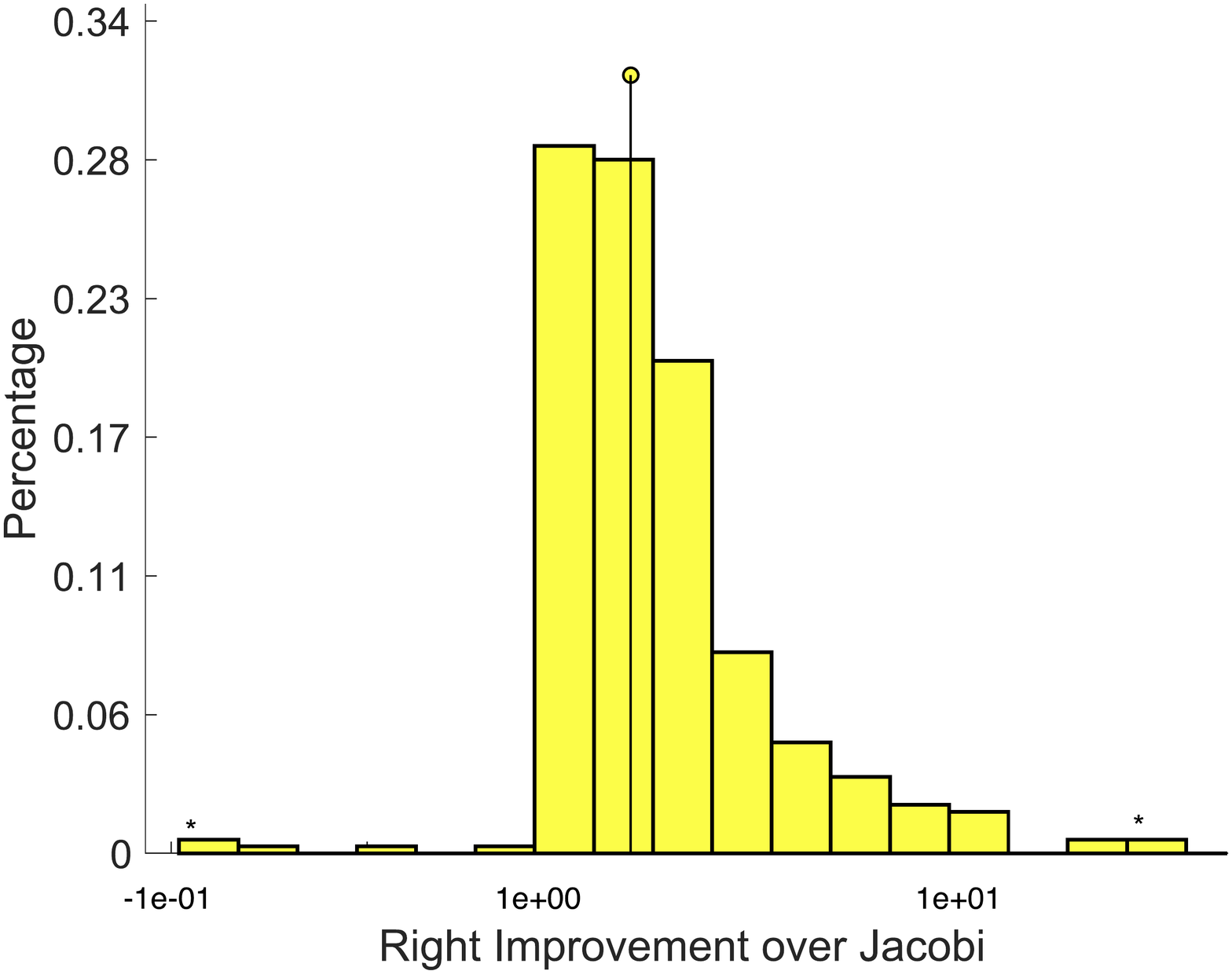}
    \includegraphics[scale=0.18]{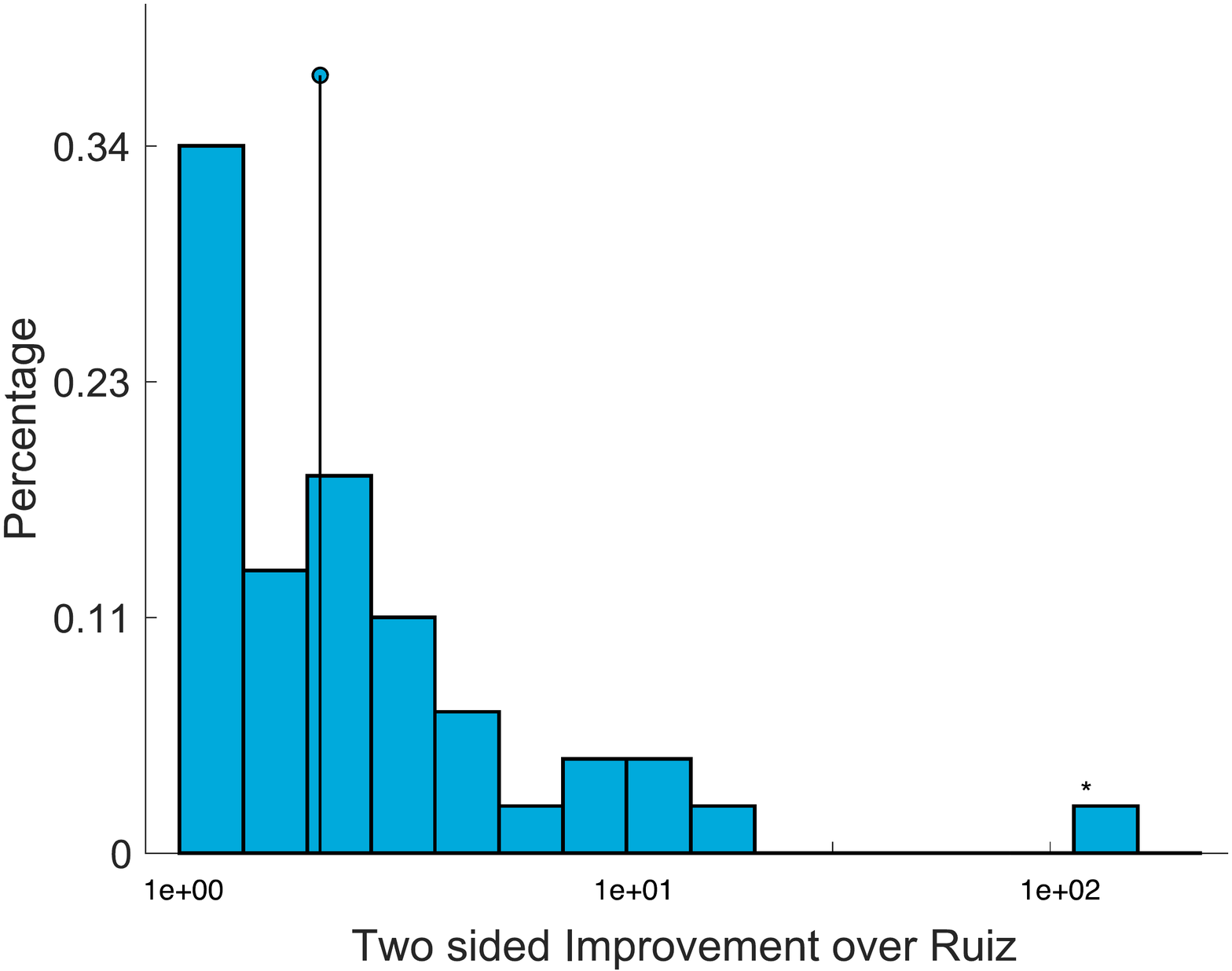}
    \caption{Comparison of log condition number improvements achieved by optimal right and two-sided preconditioners over baseline preconditioners on \texttt{SuiteSparse} matrices}
    \label{fig:comparison}
\end{figure}

As the results suggest, the dual SDP approach improves the condition number of half of the matrices by at least a factor of 2.00, and the median of the improvement is 1.49. Additionally, in Figure \ref{fig:reduction-pdf}, we plot the distribution of improvements in condition numbers for \texttt{SuiteSparse} matrices from optimal preconditioners. Moreover, compared to the baseline (Jacobi) diagonal preconditioning, the optimal preconditioning
achieves a further median improvement of 1.29. See Figure \ref{fig:comparison} for the distribution of improvements over Jacobi preconditioning. Among the \texttt{SuiteSparse} datasets, we also find a large sparse
matrix \texttt{m3plates}, whose condition number diverges after
Jacobi preconditioning, while our optimal preconditioning reduces the condition number to
close to 1.

\begin{figure}[h]
\centering
{\includegraphics[scale=0.22]{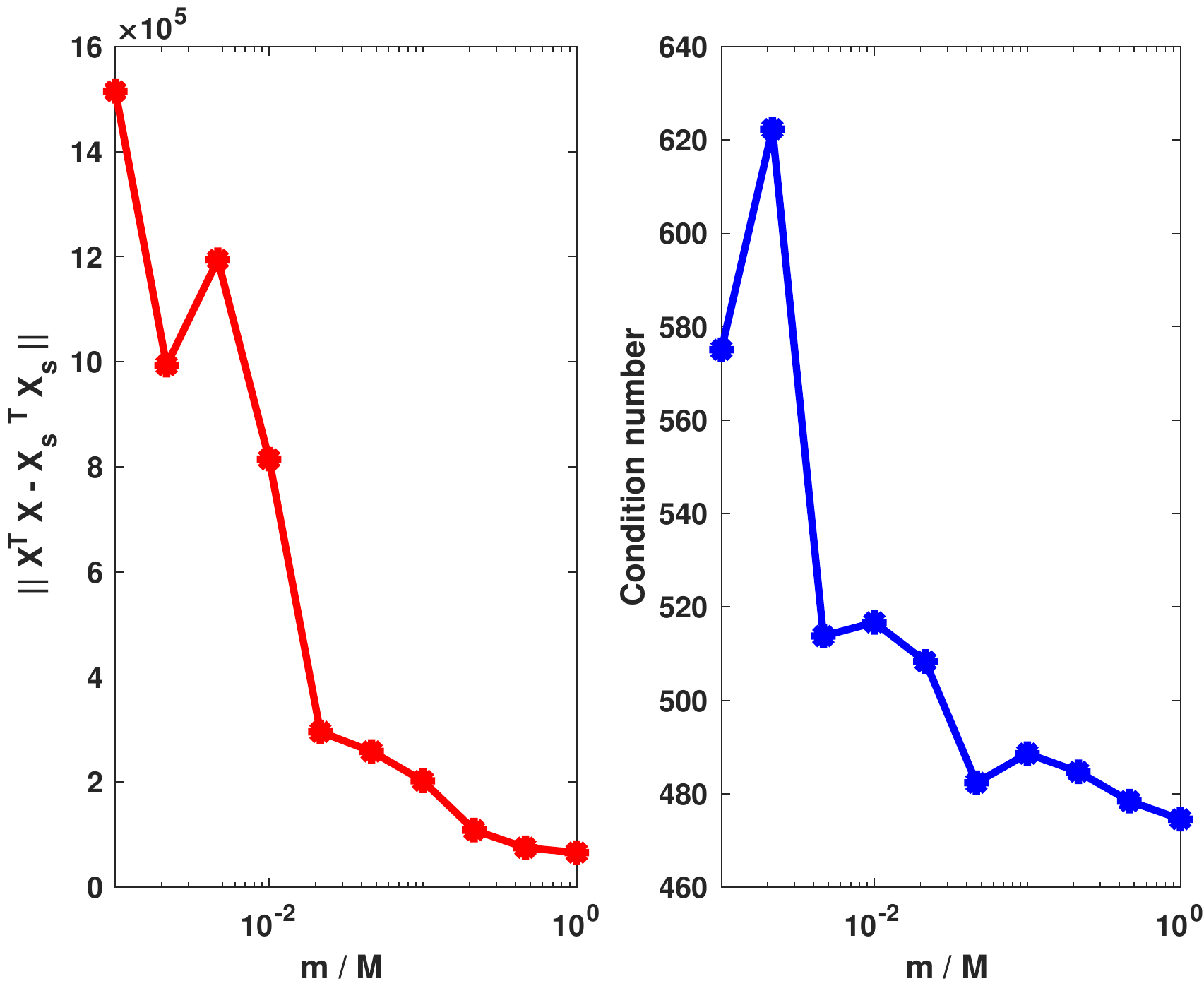}}{\includegraphics[scale=0.22]{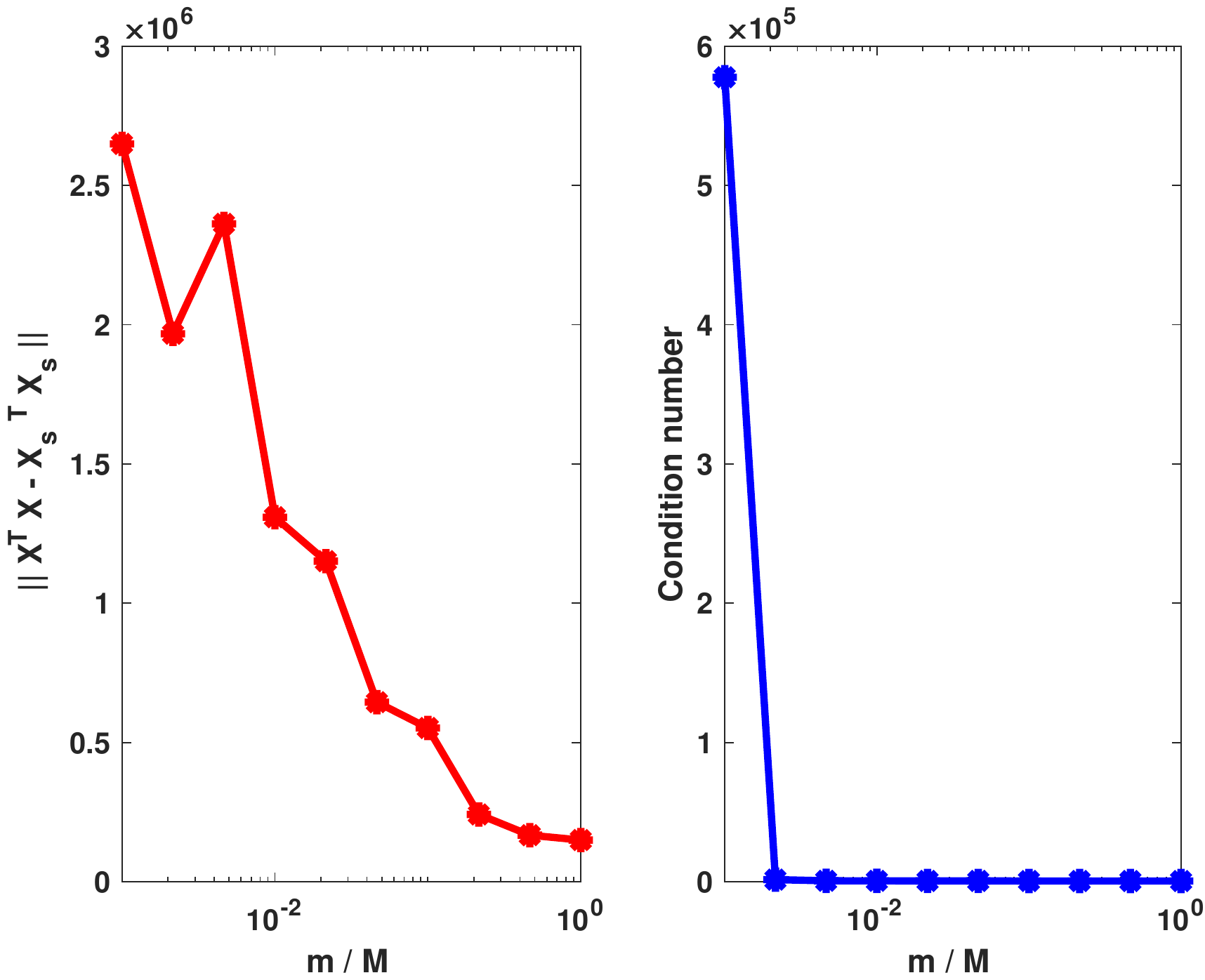}}{\includegraphics[scale=0.22]{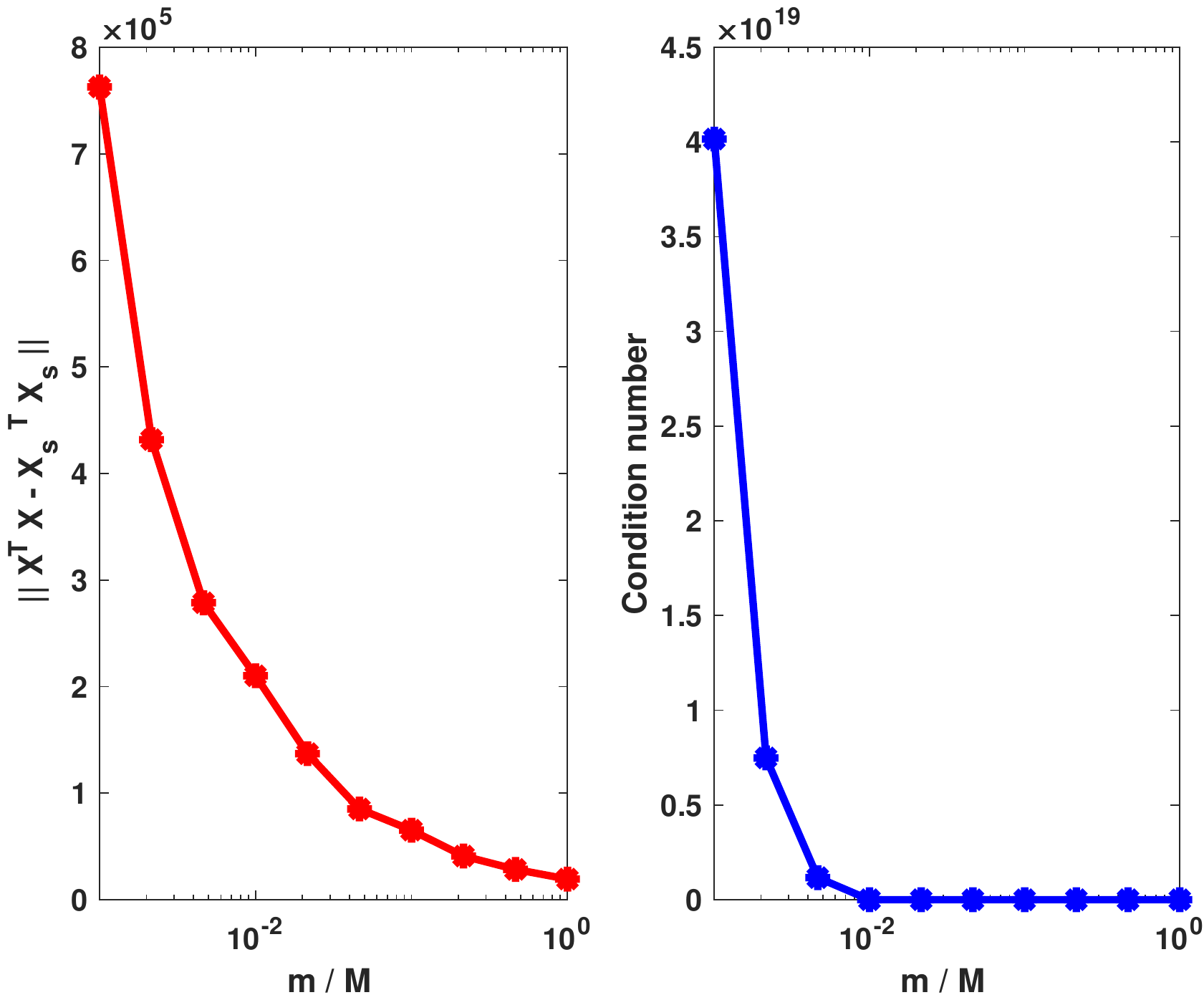}}{\includegraphics[scale=0.22]{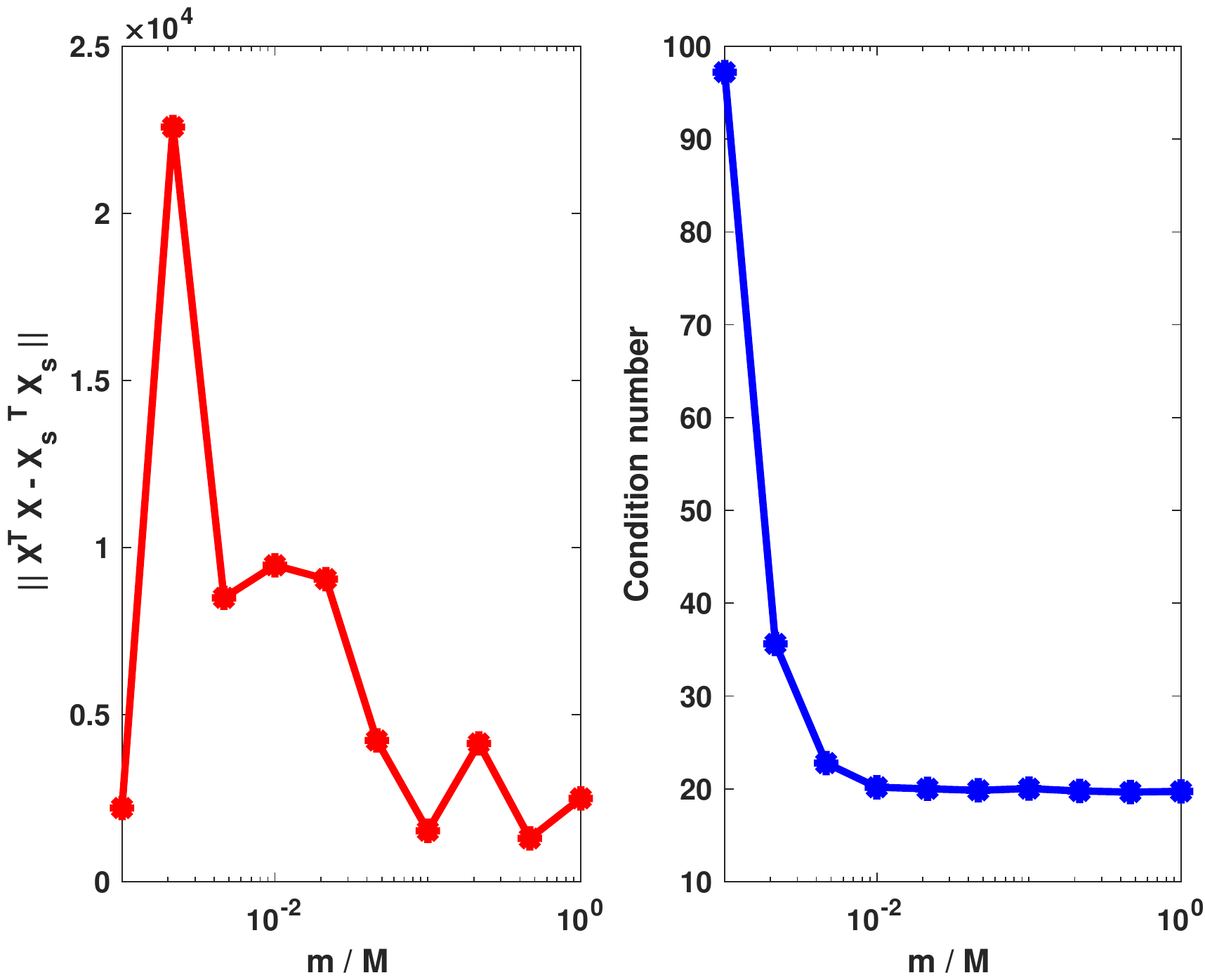}}
  \caption{Preconditioning with sampling on
  \texttt{YearPredictionMSD}, \texttt{YearPredictionMSD.t},
  \texttt{MNIST} and \texttt{credit} \ }
  \label{fig:sampling}
\end{figure}

\begin{table}[h]
\centering
  \begin{tabular}{ccccccc}
    \hline
    Mat & Size & Cbef & Caft & Cdiag & Imp. & Time\\
    \hline
    m3plates (perturbed) & 11107 & 2.651e+04 & 1.02 & $\infty$ & $>10^4$ &
    492.1\\
    \hline
  \end{tabular}
  \caption{A special example from \texttt{SuiteSparse}}
\end{table}
\paragraph{Randomized Sampling.}
In this part, we demonstrate the effectiveness of the randomization procedure on four large
real regression datasets from \texttt{LIBSVM} and \texttt{OPENML}. We summarize the results in Figure \ref{fig:sampling}.
The $x$-axis always indicates the sub-sampling ratio $\tilde{m} / m$; the red line plots $\|
A^T A - \tilde{A}^T \tilde{A} \|$ and the blue line plots the condition number
of the preconditioned matrix, using optimal right preconditioners for $A$ based on solving the dual SDP problems for the \emph{sub-sampled} matrix $\tilde A$.

\vspace{-1cm}
As the results suggest, we generally only need to randomly sample 1\% to 10\% of the rows of the original matrix $A$ in order to give
a good preconditioner for the iterative solvers. This greatly extends the
applicability of our proposed approach via solving dual SDPs, since we can efficiently find high quality right preconditioners for matrices of size 50,000 or larger (280,000 in the case of the credit dataset in \texttt{OPENL}) via random row sampling.

\paragraph{Extreme Cases} In this part, we provide random instances where the time to set up $A^TA$ is dominant and our preconditioning procedure with random row sampling can significantly reduce the condition number efficiently with less computation time.
\begin{table}[h]
\centering
  \begin{tabular}{cccccc}
    \hline
    $m$ & $n$ & Time of $A^T A$ & Time of $\tilde{A}^T \tilde{A}$ +
    computing $D_2^{-1 / 2}$  & Cbef & CAft \\
    \hline
    $10^7$ & 5 & 0.053 & 0.0004 + 0.0250 & 9.9844e+05 & 5.4020\\
    $10^7$ & 10 & 0.0999 & 0.0006 + 0.0500 & 9.9968e+05 & 1.0713\\
    $10^7$ & 20 & 0.6589 & 0.0024 + 0.4564 & 9.9889e+05 & 1.1500\\
    \hline
  \end{tabular}
  \caption{Some large random instances where randomization is efficient}
\end{table}
\subsection{Experiment Results for Left Preconditioning}
Our experiments for left preconditioning mainly focus on evaluating the reduction in condition number from left preconditioning. Since we cannot use the regularization $M+\varepsilon I$ for extremely ill-conditioned matrices as we do in right preconditioning, we present the results on 89 matrices whose associated dual SDP problems \eqref{eq:DSDP-left} are solved successfully by our solver. The summary statistics are given in Table \ref{tab:suite-left}. Detailed results can be found in the online repository.

\begin{table}[H]
\centering
  \begin{tabular}{c|c}
    \hline
    Improvement & Number\\
    \hline
    $\geq$5.00 & 31\\
    $\geq$2.00 & 55\\
    $\geq$1.25 & 81\\
    \hline
  \end{tabular}\quad\begin{tabular}{r|c}
    \hline
    Median imp. & 2.51 \\
    \hline
    Average time & 1.10\\
    \hline
  \end{tabular}
  \caption{Summary for left preconditioning on the \texttt{SuiteSparse} dataset}
  \label{tab:suite-left}
\end{table}
\vspace{-0.8cm}
Compared to right preconditioning, the dual SDPs for left preconditioning are more sensitive to the original data matrix $A$. Still, we also observe remarkable reductions in condition numbers on more than half of the tested matrices. See Figure \ref{fig:reduction-pdf} for the distribution of improvements.

\subsection{Experiment Results for Two-sided Preconditioning}
Our experiments for two-sided preconditioning compare the performance of our optimal two-sided preconditioners to the popular two-sided preconditioners based on the Ruiz scaling method \citep{ruiz2001scaling}. Since computing the two-sided preconditioner involves solving several feasibility SDPs, we restrict our attention to smaller matrices. See Table \ref{tab:two-sided} for summarized results.

\begin{table}[h]
\centering
  \begin{tabular}{c|c}
    \hline
    Improvement & Number\\
    \hline
    $\geq$5.00 & 18\\
    $\geq$2.00 & 35\\
    $\geq$1.25 & 43\\
    \hline
  \end{tabular}\quad\begin{tabular}{r|c}
    \hline
    Median imp. & 3.60 \\
    \hline
    Median imp. over Ruiz & 2.10 \\
    \hline
  \end{tabular}
  \caption{Summary for two-sided preconditioning on the \texttt{SuiteSparse} dataset}
  \label{tab:two-sided}
\end{table}
As the experiment results suggest, the optimal two-sided diagonal preconditioners work at least as well as the Ruiz scaling and, on average, an improvement of over 30\% is observed on the tested data. Moreover, the two-sided preconditioners also achieve more significant improvements in condition number compared to one-sided preconditioners (see repository for detailed comparisons). Interestingly, the Ruiz scaling does not improve much on the condition number at all for many instances. This suggests that may be room for significant improvements over current practices.

\subsection{Experiment Results for Iterative Methods}

\begin{figure}
    \centering
    \includegraphics[scale=0.3]{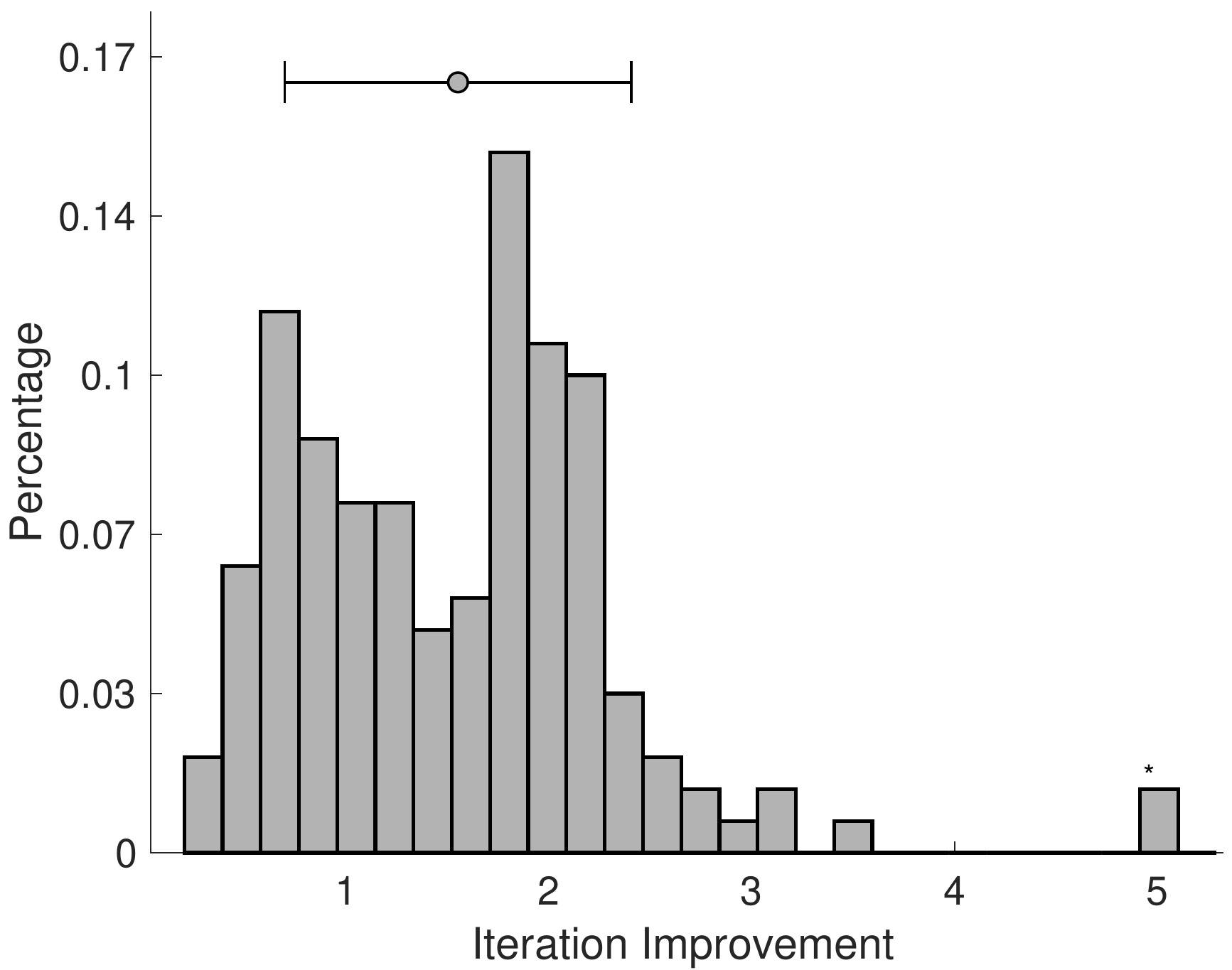}
    \caption{Improvements in iteration numbers of the preconditioned conjugate gradient method on the testing collection using optimal right diagonal preconditioners vs. Jacobi preconditioners}
    \label{fig:iterative-comparison}
\end{figure}

In this section, we evaluate the practical efficiency of our optimal preconditioner at speeding up iterative methods. We focus on the right-preconditioners whose improvements are greater than 2 and use the number of preconditioned conjugate gradient iterations as the evaluation metric.
The testing matrices are selected such that $\kappa (D_2^{- 1 / 2} M D_2^{- 1 / 2}) \ll \kappa (D_M^{- 1 / 2} M D_M^{- 1 / 2})$. We end up with a collection of 148 matrices. In Table \ref{tab:iterative-comparison}, we report the number of instances for which the optimal and
the baseline Jacobi preconditioners perform better.
\begin{table}[h]
\centering
\begin{tabular}{r|c}
    \hline
    Preconditioner & \# of instances\\
    \hline
    Jacobi & 46\\
    \hline
    Optimal & 102\\
    \hline
    Total & 148\\
    \hline
  \end{tabular}
  \caption{Instances where each preconditioning technique yields fastest convergence}
  \label{tab:iterative-comparison}
\end{table}
Figure \ref{fig:iterative-comparison} illustrates the relative improvements in terms of number of iterations using optimal diagonal preconditioners over Jacobi preconditioners.


\subsection{Comparing Iterative and Direct Methods}
In this section, we further test the applicability of optimal preconditioners to iterative methods as a potential substitute for sparse direct methods, especially if we need to solve linear systems with identical LHS matrix and varying RHSs. We simulate the scenario, given a matrix $A$, as follows: 
\begin{enumerate}
    \item Compute the optimal diagonal left/right preconditioner and explicitly evaluate the pre-conditioned matrix.
    \item Compute the symbolic ordering $P$ and factorize $P^TA^TAP= LL^T$.
\end{enumerate}
Since the time of the above two operations can be amortized by the large number of solves, we exclude them from our time statistics. Also, we explicitly compute the preconditioned matrix since 1) diagonal preconditioners have no impact on the sparsity pattern of the matrix, and 2) computational time for the preconditioning step is amortized.

Then 12800 linear systems with identical LHSs are solved and the time statistics for both sparse and direct methods are collected. In our setup, we adopt the native routine from \texttt{MATLAB} to perform forward/backward substitution and the \texttt{Intel MKL Reverse Communication Interface} to implement the CG solver. We found that in 35 instances, right optimally preconditioned iterative solvers are faster than direct solvers, and in 14 instances left optimally-conditioned preconditioners are faster than direct solvers. See the Github repository for detailed experimental results. Even though our primary focus is on maximal condition number reduction with diagonal preconditioners, we have demonstrated that they can also speed up iterative methods significantly and, in some cases, outperform sparse direct methods on instances from \texttt{SuiteSparse}.

\subsection{Scalability of the Dual SDP Approach}

To conclude this section, we make a few remarks on the scalability of our dual SDP
solver in solving optimal diagonal preconditioning problems for $A \in \mathbb{R}^{m \times n}$ 
as a reference for practical
use.
\begin{enumerate}[leftmargin=15pt]
  \item The complexity (or CPU time) of the three SDP formulations ranks as
  Two-sided $\gg$ Left $\geq$ Right, where $=$ holds if $m \approx n$. Correspondingly, the reduction in condition number ranks as Two-sided $>$
  Right $\approx$ Left. 
  
  Although two-sided scaling is most powerful and results in the greatest
  reduction in condition numbers, its high computation cost restricts its use
  to matrices with $\max \{ m, n \} \leq 300$. But as our computational
  experience suggests, left or right preconditioners suffice in most cases.
  
  \item Left and Right preconditioners benefit from sparsity.
  
  If $A$ is extremely sparse, SDP solvers could exploit this to accelerate the solution of the dual SDPs \eqref{eq:DSDP-right} and \eqref{eq:DSDP-left}. In our experiments, sparsity makes it possible to compute preconditioners of
  a $10000 \times 10000$ matrix very efficiently. In general cases, solvers of the dual SDP \eqref{eq:DSDP-left} for left
  preconditioners can handle matrices $A \in \mathbb{R}^{m \times n}$ satisfying $\max \{ m, n \} \leq 2000$ and solvers of the dual SDP \eqref{eq:DSDP-right} for
  right preconditioners can handle cases where $n \leq 5000$.
  
  \item Right preconditioners benefit from statistical properties.
  
  If the rows of $A$ are generated from some statistical distribution, as is the case for data matrices in practice, then right preconditioners can
  handle very tall matrices with large $m$ quite efficiently using the aforementioned random sampling
  procedure. Especially when $m \gg n =\mathcal{O} (1)$, a right
  preconditioner with high quality can be found even faster than computing $A^T A$.
\end{enumerate}

\section{Concluding Remarks}
In this paper, we study the problem of optimal diagonal preconditioning for any matrix $A$ with full rank. More precisely, the goal is to achieve maximal condition number reduction of $A$ by scaling its rows and/or columns. This problem is relevant to the optimization and OR community both from a scientific standpoint and from a practical standpoint. First, as existing popular diagonal preconditioners often do not have guarantees on condition number reductions, it is important to know to what extent we may improve upon existing heuristics-based diagonal preconditioners. Second, since condition numbers often affect the convergence speed of algorithms, we are interested in efficient methods to find near optimal diagonal preconditioners for large systems, so that they could be applied to speed up the convergence. Starting from a quasi-convex reformulation of the general two-sided optimal diagonal preconditioning problem, we provide bisection and interior point algorithms based on potential reduction, which are also applicable to other quasi-convex problems with PSD matrix variables. Specializing to one-sided optimal diagonal preconditioning problems, we demonstrate they are equivalent to standard dual SDP problems with sparse and low rank structures. We then develop efficient customized solvers incorporating randomization, which are able to handle problems of size as large as 200,000 efficiently. Based on our extensive experiments, we demonstrate that there is indeed a large room of improvement over existing diagonal preconditioners in terms of condition number reduction. Moreover, on large real datasets, our optimal diagonal preconditioners with customized solvers are shown to improve over heuristic preconditioners at speeding up iterative methods and can also outperform sparse direct methods when the problem requires solving linear systems with identical LHS and multiple RHSs. Although the computational overhead of finding optimal diagonal preconditioners is a potential concern, many modern optimization settings can benefit from optimal diagonal preconditioners of smaller subsystems, such as in ADMM or Federated Learning, or a system that needs to be repeatedly solved, such as preconditioned CG in the Newton step of interior point methods or solving multi-scale numerical PDEs. We therefore see our results as both theoretically and scientifically interesting, and also practically relevant to problems that may benefit from better preconditioning. 

\ACKNOWLEDGMENT{The authors gratefully acknowledge helpful discussions with Professor Stephen Boyd, Wenbo Gao, and participants of the Optimization Group Meeting organized by Professor Yinyu Ye.}


\bibliographystyle{informs2014} 
\bibliography{preconditioning_references} 

\begin{thebibliography}{49}
\providecommand{\natexlab}[1]{#1}
\providecommand{\url}[1]{\texttt{#1}}
\providecommand{\urlprefix}{URL }

\bibitem[{Andersen et~al.(2011)Andersen, Dahl, Liu, Vandenberghe, Sra, Nowozin,
  \protect\BIBand{} Wright}]{andersen2011interior}
Andersen M, Dahl J, Liu Z, Vandenberghe L, Sra S, Nowozin S, Wright S (2011)
  Interior-point methods for large-scale cone programming. \emph{Optimization
  for machine learning} 5583.

\bibitem[{Barrett et~al.(1994)Barrett, Berry, Chan, Demmel, Donato, Dongarra,
  Eijkhout, Pozo, Romine, \protect\BIBand{} Van~der
  Vorst}]{barrett1994templates}
Barrett R, Berry M, Chan TF, Demmel J, Donato J, Dongarra J, Eijkhout V, Pozo
  R, Romine C, Van~der Vorst H (1994) \emph{Templates for the solution of
  linear systems: building blocks for iterative methods} (SIAM).

\bibitem[{Beck \protect\BIBand{} Teboulle(2003)}]{beck2003mirror}
Beck A, Teboulle M (2003) Mirror descent and nonlinear projected subgradient
  methods for convex optimization. \emph{Operations Research Letters}
  31(3):167--175.

\bibitem[{Benson(1982)}]{benson1982iterative}
Benson MW (1982) Iterative solution of large sparse linear systems arising in
  certain multidimensional approximation problems. \emph{Util. Math.}
  22:127--140.

\bibitem[{Benson \protect\BIBand{} Ye(2008)}]{benson2008algorithm}
Benson SJ, Ye Y (2008) Algorithm 875: Dsdp5—software for semidefinite
  programming. \emph{ACM Transactions on Mathematical Software (TOMS)}
  34(3):1--20.

\bibitem[{Benzi(2002)}]{benzi2002preconditioning}
Benzi M (2002) Preconditioning techniques for large linear systems: a survey.
  \emph{Journal of computational Physics} 182(2):418--477.

\bibitem[{Benzi \protect\BIBand{} Tuma(1999)}]{benzi1999comparative}
Benzi M, Tuma M (1999) A comparative study of sparse approximate inverse
  preconditioners. \emph{Applied Numerical Mathematics} 30(2-3):305--340.

\bibitem[{Boyd et~al.(2011)Boyd, Parikh, Chu, Peleato, Eckstein
  et~al.}]{boyd2011distributed}
Boyd S, Parikh N, Chu E, Peleato B, Eckstein J, et~al. (2011) Distributed
  optimization and statistical learning via the alternating direction method of
  multipliers. \emph{Foundations and Trends{\textregistered} in Machine
  learning} 3(1):1--122.

\bibitem[{Bradley(2010)}]{bradley2010algorithms}
Bradley AM (2010) Algorithms for the equilibration of matrices and their
  application to limited-memory quasi-newton methods. Technical report,
  STANFORD UNIV CA.

\bibitem[{Broyden(1970)}]{broyden1970convergence}
Broyden CG (1970) The convergence of a class of double-rank minimization
  algorithms 1. general considerations. \emph{IMA Journal of Applied
  Mathematics} 6(1):76--90.

\bibitem[{Burer et~al.(2002)Burer, Monteiro, \protect\BIBand{}
  Zhang}]{burer2002solving}
Burer S, Monteiro RD, Zhang Y (2002) Solving a class of semidefinite programs
  via nonlinear programming. \emph{Mathematical Programming} 93(1):97--122.

\bibitem[{Cesari(1937)}]{cesari1937sulla}
Cesari L (1937) \emph{Sulla risoluzione dei sistemi di equazioni lineari per
  approssimazioni successive} (G. Bardi).

\bibitem[{Chang \protect\BIBand{} Lin(2011)}]{chang2011libsvm}
Chang CC, Lin CJ (2011) Libsvm: a library for support vector machines.
  \emph{ACM transactions on intelligent systems and technology (TIST)}
  2(3):1--27.

\bibitem[{Concus et~al.(1985)Concus, Golub, \protect\BIBand{}
  Meurant}]{concus1985block}
Concus P, Golub GH, Meurant G (1985) Block preconditioning for the conjugate
  gradient method. \emph{SIAM Journal on Scientific and Statistical Computing}
  6(1):220--252.

\bibitem[{Dubois et~al.(1979)Dubois, Greenbaum, \protect\BIBand{}
  Rodrigue}]{dubois1979approximating}
Dubois PF, Greenbaum A, Rodrigue GH (1979) Approximating the inverse of a
  matrix for use in iterative algorithms on vector processors. \emph{Computing}
  22(3):257--268.

\bibitem[{Gao et~al.(2022)Gao, Ge, \protect\BIBand{} Ye}]{gao2022hdsdp}
Gao W, Ge D, Ye Y (2022) Hdsdp: Software for semidefinite programming.
  \emph{arXiv preprint arXiv:2207.13862} .

\bibitem[{Gao \protect\BIBand{} Goldfarb(2019)}]{gao2019quasi}
Gao W, Goldfarb D (2019) Quasi-newton methods: superlinear convergence without
  line searches for self-concordant functions. \emph{Optimization Methods and
  Software} 34(1):194--217.

\bibitem[{Goldfarb(1970)}]{goldfarb1970family}
Goldfarb D (1970) A family of variable-metric methods derived by variational
  means. \emph{Mathematics of computation} 24(109):23--26.

\bibitem[{Hestenes \protect\BIBand{} Stiefel(1952)}]{hestenes1952methods}
Hestenes MR, Stiefel E (1952) Methods of conjugate gradients for solving.
  \emph{Journal of research of the National Bureau of Standards} 49(6):409.

\bibitem[{Jacobi(1845)}]{jacobi1845ueber}
Jacobi CG (1845) Ueber eine neue aufl{\"o}sungsart der bei der methode der
  kleinsten quadrate vorkommenden line{\"a}ren gleichungen. \emph{Astronomische
  Nachrichten} 22(20):297--306.

\bibitem[{Jambulapati et~al.(2020)Jambulapati, Li, Musco, Sidford,
  \protect\BIBand{} Tian}]{jambulapati2020fast}
Jambulapati A, Li J, Musco C, Sidford A, Tian K (2020) Fast and near-optimal
  diagonal preconditioning. \emph{arXiv preprint arXiv:2008.01722} .

\bibitem[{Johnson et~al.(1983)Johnson, Micchelli, \protect\BIBand{}
  Paul}]{johnson1983polynomial}
Johnson OG, Micchelli CA, Paul G (1983) Polynomial preconditioners for
  conjugate gradient calculations. \emph{SIAM Journal on Numerical Analysis}
  20(2):362--376.

\bibitem[{Kairouz et~al.(2021)Kairouz, McMahan, Avent, Bellet, Bennis, Bhagoji,
  Bonawitz, Charles, Cormode, Cummings et~al.}]{kairouz2021advances}
Kairouz P, McMahan HB, Avent B, Bellet A, Bennis M, Bhagoji AN, Bonawitz K,
  Charles Z, Cormode G, Cummings R, et~al. (2021) Advances and open problems in
  federated learning. \emph{Foundations and Trends{\textregistered} in Machine
  Learning} 14(1--2):1--210.

\bibitem[{Kaszkurewicz et~al.(1995)Kaszkurewicz, Bhaya, \protect\BIBand{}
  Ramos}]{kaszkurewicz1995control}
Kaszkurewicz E, Bhaya A, Ramos P (1995) A control-theoretic view of diagonal
  preconditioners. \emph{International journal of systems science}
  26(9):1659--1672.

\bibitem[{Kolodziej et~al.(2019)Kolodziej, Aznaveh, Bullock, David, Davis,
  Henderson, Hu, \protect\BIBand{} Sandstrom}]{kolodziej2019suitesparse}
Kolodziej SP, Aznaveh M, Bullock M, David J, Davis TA, Henderson M, Hu Y,
  Sandstrom R (2019) The suitesparse matrix collection website interface.
  \emph{Journal of Open Source Software} 4(35):1244.

\bibitem[{Kone{\v{c}}n{\`y} et~al.(2016)Kone{\v{c}}n{\`y}, McMahan, Yu,
  Richt{\'a}rik, Suresh, \protect\BIBand{} Bacon}]{konevcny2016federated}
Kone{\v{c}}n{\`y} J, McMahan HB, Yu FX, Richt{\'a}rik P, Suresh AT, Bacon D
  (2016) Federated learning: Strategies for improving communication efficiency.
  \emph{arXiv preprint arXiv:1610.05492} .

\bibitem[{Lanczos(1950)}]{lanczos1950iteration}
Lanczos C (1950) An iteration method for the solution of the eigenvalue problem
  of linear differential and integral operators .

\bibitem[{Liesen \protect\BIBand{} Strakos(2013)}]{liesen2013krylov}
Liesen J, Strakos Z (2013) \emph{Krylov subspace methods: principles and
  analysis} (Oxford University Press).

\bibitem[{Lin et~al.(2018)Lin, Ma, Ye, \protect\BIBand{} Zhang}]{lin2018admm}
Lin T, Ma S, Ye Y, Zhang S (2018) An admm-based interior-point method for
  large-scale linear programming. \emph{arXiv preprint arXiv:1805.12344} .

\bibitem[{Luenberger et~al.(1984)Luenberger, Ye et~al.}]{luenberger1984linear}
Luenberger DG, Ye Y, et~al. (1984) \emph{Linear and nonlinear programming},
  volume~2 (Springer).

\bibitem[{Luo et~al.(1998)Luo, Sturm, \protect\BIBand{}
  Zhang}]{luo1998superlinear}
Luo ZQ, Sturm JF, Zhang S (1998) Superlinear convergence of a symmetric
  primal-dual path following algorithm for semidefinite programming. \emph{SIAM
  Journal on Optimization} 8(1):59--81.

\bibitem[{Mizuno(1992)}]{mizuno1992new}
Mizuno S (1992) A new polynomial time method for a linear complementarity
  problem. \emph{Mathematical Programming} 56(1-3):31--43.

\bibitem[{Mizuno et~al.(1993)Mizuno, Todd, \protect\BIBand{}
  Ye}]{mizuno1993adaptive}
Mizuno S, Todd MJ, Ye Y (1993) On adaptive-step primal-dual interior-point
  algorithms for linear programming. \emph{Mathematics of Operations research}
  18(4):964--981.

\bibitem[{Nesterov \protect\BIBand{} Todd(1997)}]{nesterov1997self}
Nesterov YE, Todd MJ (1997) Self-scaled barriers and interior-point methods for
  convex programming. \emph{Mathematics of Operations research} 22(1):1--42.

\bibitem[{Nocedal \protect\BIBand{} Wright(2006)}]{nocedal2006numerical}
Nocedal J, Wright S (2006) \emph{Numerical optimization} (Springer Science \&
  Business Media).

\bibitem[{Rijn et~al.(2013)Rijn, Bischl, Torgo, Gao, Umaashankar, Fischer,
  Winter, Wiswedel, Berthold, \protect\BIBand{} Vanschoren}]{rijn2013openml}
Rijn JNv, Bischl B, Torgo L, Gao B, Umaashankar V, Fischer S, Winter P,
  Wiswedel B, Berthold MR, Vanschoren J (2013) Openml: A collaborative science
  platform. \emph{Joint european conference on machine learning and knowledge
  discovery in databases}, 645--649 (Springer).

\bibitem[{Ruiz(2001)}]{ruiz2001scaling}
Ruiz D (2001) A scaling algorithm to equilibrate both rows and columns norms in
  matrices. Technical report, CM-P00040415.

\bibitem[{Saad(2003)}]{saad2003iterative}
Saad Y (2003) \emph{Iterative methods for sparse linear systems}, volume~82
  (siam).

\bibitem[{Sinkhorn(1964)}]{sinkhorn1964relationship}
Sinkhorn R (1964) A relationship between arbitrary positive matrices and doubly
  stochastic matrices. \emph{The annals of mathematical statistics}
  35(2):876--879.

\bibitem[{Sturm \protect\BIBand{} Zhang(1999)}]{sturm1999symmetric}
Sturm JF, Zhang S (1999) Symmetric primal-dual path-following algorithms for
  semidefinite programming. \emph{Applied Numerical Mathematics}
  29(3):301--315.

\bibitem[{Takapoui \protect\BIBand{}
  Javadi(2016)}]{takapoui2016preconditioning}
Takapoui R, Javadi H (2016) Preconditioning via diagonal scaling. \emph{arXiv
  preprint arXiv:1610.03871} .

\bibitem[{Todd et~al.(1998)Todd, Toh, \protect\BIBand{}
  T{\"u}t{\"u}nc{\"u}}]{todd1998nesterov}
Todd MJ, Toh KC, T{\"u}t{\"u}nc{\"u} RH (1998) On the nesterov--todd direction
  in semidefinite programming. \emph{SIAM Journal on Optimization}
  8(3):769--796.

\bibitem[{Turing(1948)}]{turing1948rounding}
Turing AM (1948) Rounding-off errors in matrix processes. \emph{The Quarterly
  Journal of Mechanics and Applied Mathematics} 1(1):287--308.

\bibitem[{Von~Neumann(1937)}]{von1937uber}
Von~Neumann J (1937) Uber ein okonomsiches gleichungssystem und eine
  verallgemeinering des browerschen fixpunktsatzes. \emph{Erge. Math. Kolloq.},
  volume~8, 73--83.

\bibitem[{Von~Neumann(1945)}]{neumann1945model}
Von~Neumann J (1945) A model of general economic equilibrium. \emph{The Review
  of Economic Studies} 13(1):1--9.

\bibitem[{Ye(1992)}]{ye1992potential}
Ye Y (1992) A potential reduction algorithm allowing column generation.
  \emph{SIAM Journal on Optimization} 2(1):7--20.

\bibitem[{Ye(1995)}]{ye1995neumann}
Ye Y (1995) On the von neumann economic growth problem. \emph{Mathematics of
  operations research} 20(3):617--633.

\bibitem[{Young(1954)}]{young1954iterative}
Young D (1954) Iterative methods for solving partial difference equations of
  elliptic type. \emph{Transactions of the American Mathematical Society}
  76(1):92--111.

\bibitem[{Young(2014)}]{young2014iterative}
Young DM (2014) \emph{Iterative solution of large linear systems} (Elsevier).

\end{thebibliography}



%
%
%
 \begin{APPENDICES}
 \newpage
\section{Basic Facts about the Condition Number}
\label{sec:appendix_basic}

Let $\|\cdot\|$ be any matrix norm on the space of matrices $\mathbb{R}^{m\times n}$.
Then the condition number associated with this norm is defined as
\begin{align*}
\kappa(A): & =\|A\|\cdot\|A^{-1}\|
\end{align*}
where $A^{-1}$ is the pseudo-inverse of a non-singular $A$. 

When $\|\cdot\|$ is the $2$-norm defined as 
\begin{align*}
\|A\| & =\sup_{\|v\|_{2}=1}\|Av\|_{2}
\end{align*}
the condition number can be alternatively written 
\begin{align*}
\kappa(A) & =\frac{\max_{i}\sigma_{i}(A)}{\min_{i}\sigma_{i}(A)}
\end{align*}
where $\{\sigma_{i}\}_{i=1}^{m\land n}$ are the non-zero singular
values of $A$. Moreover, if $A$ is symmetric positive definite,
then we can replace singular values with eigenvalues in the definition. 

The condition number is a measure of how close a matrix is to being
singular. Alternatively, it also describes the amount of distortion
that the matrix transformation causes when applied to the unit sphere. 

Because matrix norms are equivalent, in the sense that for any two
norms $\|\cdot\|_{a}$ and $\|\cdot\|_{b}$, there exist constants
$c,C$ such that
\begin{align*}
c\|A\|_{a} & \leq\|A\|_{b}\leq C\|A\|_{a}
\end{align*}
we see that the condition numbers defined by the two norms are also
equivalent, i.e. 
\begin{align*}
c^{2}\kappa_{a}(A)\leq & \kappa_{b}(A)\leq C^{2}\kappa_{a}(A)
\end{align*}

Condition number is invariant under scaling by constants, and from
definition, $\kappa(A)=\kappa(A^{-1})$. Since $\|AB\|\leq\|A\|\|B\|$
for matrix norms, $\kappa(AB)\leq\kappa(A)\cdot\kappa(B)$, and this
also implies $\kappa(A)\geq1$ for any matrix $A$.

Lastly, for condition numbers defined with matrix $p$-norms, we have
\begin{align*}
\kappa(A^{T}A) & =\kappa(A^{T})\cdot\kappa(A)\\
& =\kappa(A)^{2}
\end{align*}
since 
\begin{align*}
\max_{\|v\|=1}\|A^{T}Av\| & =\max_{\|u\|=1,\|v\|=1}\langle u,A^{T}Av\rangle\\
& \geq\max_{\|v\|=1}\langle Av,Av\rangle=\|A\|^{2}
\end{align*}
and so in fact $\|A^{T}A\|=\|A\|^{2}$.

\section{Details on the Interior Point Algorithm in Section \ref{sec:optimal}}
\label{sec:appendix_optimal}

In this section, we provide details on the potential reduction algorithms discussed in Section \ref{sec:potential-reduction}. Recall the potential function associated with the two-sided optimal diagonal preconditioning problem \eqref{eq:SDP-two-sided-upper-bound}
\begin{align*}
P(\kappa):=\max_{D_{1},D_{2}\in \Gamma(\kappa)} & \log\det(A^{T}D_{1}A-D_{2})+\log\det(\kappa D_{2}-A^{T}D_{1}A)+\log\det(D_{1}-I_{m})+\log\det(\hat \kappa I_{m}-D_1)
\end{align*}
where $\Gamma(k)$ is the feasible region given by 
\begin{align*}
\Gamma(\kappa): & =\{D_{1},D_{2}\text{ diagonal:}\text{ }A^{T}D_{1}A\succeq D_{2};\text{ }\kappa D_{2}\succeq A^{T}D_{1}A;\text{ }\hat \kappa I_m \succeq D_{1}\succeq I_{m}\}
\end{align*}

First, note that if the feasible region $\Gamma(\kappa^{0})$ has a non-empty
interior for some $\kappa^{0}<\infty$, then the feasible region $\Gamma(\kappa)$
has a non-empty interior for all $\kappa\geq\kappa^{0}$. Since the feasible region is bounded,
the potential function is well-defined for all $\kappa> \kappa^\ast$, which is the optimal condition number achieved by diagonal preconditioners, and also the minimum value $\kappa^{\ast}\geq1$
with a non-empty feasible region. Note, however, that the feasible region $\Gamma(\kappa^{\ast})$
has no interior, because it must be the case that $\lambda_{\min}(\kappa^{\ast}D_{2}-A^{T}D_{1}A)=0$
for any $D_{1},D_{2}$ in $\Gamma(\kappa^{\ast})$. Otherwise, $\kappa$
can be decreased without violating the constraints defining the feasible
region, yielding a non-empty feasible region $\Gamma(\kappa')$ for
some $\kappa'<\kappa^{\ast}$. Lastly, if $(D_1,D_2) \in \Gamma(\kappa)$, then so is $(cD_1,cD_2)$ for any $c>0$ such that $\overline \kappa I_m \succeq cD_{1}\succeq I_{m}$, but the analytic center, which is the maximizer in the definition of $P(\kappa)$, is unique. 

Next, we show that the potential function $P(\kappa)$ is monotone increasing in $\kappa$,
and is bounded below for any $\kappa$
large. This implies that $P(\kappa)$ is a reasonable measure of the progress of an algorithm that aims to reduce the condition number:
to obtain an approximate solution to the optimal preconditioning problem, it suffices to reduce the potential function by a constant
amount at every step and terminate when the potential function is
reduced below a certain threshold determined by accuracy $\epsilon$. The resulting analytic centers are then approximately optimal preconditioners. These properties of $P(\kappa)$ are summarized in the following two results. All omitted proofs can be found in Appendix \ref{sec:proofs}.
\begin{proposition}\label{prop:monotonicity}
	If $\kappa^0>\kappa^1>\kappa^{\ast}$, then $P(\kappa^0)>P(\kappa^{1})$. 
\end{proposition}
\begin{proposition}\label{prop:lower bound}
	Let $t$ satisfy $\kappa^{\ast}\leq2^{t}$, then for any $\overline{\kappa}>\kappa^{\ast}+2^{-t}$,
	we have
	\begin{align*}
	P(\overline{\kappa}) & >-Cmt
	\end{align*}
	for some universal constant $C$ that only depends on the matrix $A\in \mathbb{R}^{m\times n}$. 
\end{proposition}
If $P(\kappa)<-Cmt$, then we must
have $\kappa\leq\kappa^{\ast}+2^{-t}$. Therefore, if we can find an analytic center 
$(D_1,D_2)$ whose associated potential is below the threshold of $O(m\log\epsilon)$,
we have found a pair of preconditioners that yields a condition number
that is within $\epsilon$ of the optimal condition number. Our goal
is then to design an algorithm that achieves constant reduction of
the potential at every iteration at a low computation cost, which would guarantee an iteration
complexity of $O(\log\frac{1}{\epsilon})$ with a low computation complexity.  

In the next two subsections, we develop potential reduction algorithms based on this idea, using exact and approximate analytic centers, respectively. 

\subsection{Potential Reduction with Exact Analytic Centers}

Let $(D_1(\kappa^0),D_2(\kappa^0)),(D_1(\kappa^1),D_2(\kappa^1))$ be the analytic centers of $\Gamma^{0}=\Gamma(\kappa^{0})$
and $\Gamma^{1}=\Gamma(\kappa^{1})$ for $\kappa^0,\kappa^1>\kappa^\ast$. With a suitable choice of $\Delta\kappa=\kappa^{0}-\kappa^{1}>0$,
we show that 
\begin{align*}
P(\kappa^{1}) & \leq P(\kappa^{0})-\Omega(1)
\end{align*}
so that we can achieve constant reductions of the potential function
by updating from $(D_1(\kappa^0),D_2(\kappa^0))$ to $(D_1(\kappa^1),D_2(\kappa^1))$ using Newton iterations with a quadratic
convergence rate. The exact form of the Newton iteration requires the Nesterov-Todd direction \citep{todd1998nesterov}, and is therefore deferred
to the next subsection, where we propose a practical algorithm
using approximate analytic centers and only one Newton update at each iteration.
\begin{theorem}
	\label{thm:exact center}Let $\overline{\kappa}$ be an upper bound of $\kappa^{\ast}$. Let $(D_1(\kappa^0),D_2(\kappa^0))$ and $(D_1(\kappa^1),D_2(\kappa^1))$ be the analytic centers of $\Gamma^{0}=\Gamma(\kappa^{0})$
	and $\Gamma^{1}=\Gamma(\kappa^{1})$ for $\kappa^{0},\kappa^{1}\in(\kappa^{\ast},\overline{\kappa})$.
	If $\Delta\kappa=\kappa^{0}-\kappa^{1}>0$ satisfies
	\begin{align*}
\Delta\kappa & =\frac{\beta}{Tr(D_{2}(\kappa^{0})\cdot(\kappa^{0}D_{2}(\kappa^{0})-A^{T}D_{1}(\kappa^{0})A)^{-1})}
	\end{align*}
	for some sufficiently small $\beta\in(0,1)$ that only depends on
	$\overline{\kappa},\hat{\kappa},A$, then 
	\begin{align*}
	P(\kappa^{1}) \leq P(\kappa^{0})-\beta. 
	\end{align*}
\end{theorem}

One can also show that the analytic centers $(D_1(\kappa^0),D_2(\kappa^0))$ and $(D_1(\kappa^1),D_2(\kappa^1))$ are close to each other, although this property is not needed at present. Theorem \ref{thm:exact center} suggests the following conceptual algorithm that relies on moving along the central path of analytic centers by decreasing $\kappa$ by $\Delta \kappa$ at each step.
\begin{algorithm}[Potential Reduction with Exact Analytic Centers]
	\label{alg:exact}
\begin{enumerate}
	\normalfont 
	    \item Let $(D_1(\kappa^0),D_2(\kappa^0))$ be the analytic center of $\Gamma(\kappa^{0})$ with some
	$\kappa^{0}>\kappa^{\ast}$. Set $k=0$. 
	
	\item  With a small fixed constant $\beta$, set	\begin{align*}
\Delta\kappa^k & =\frac{\beta}{Tr(D_{2}(\kappa^{k})\cdot(\kappa^{k}D_{2}(\kappa^{k})-A^{T}D_{1}(\kappa^{k})A)^{-1})},
	\end{align*}
and let $\kappa^{k+1}=\kappa^{k}-\Delta\kappa^{k}$.
	Then apply Newton iterations described in Proposition \ref{prop:Newton step} to generate
	the exact analytic center of $\Gamma(\kappa^{k+1})$ from the analytic center of $\Gamma(\kappa^k)$. 
	
	\item Compute $P(\kappa^{k+1})$. If $P(\kappa^{k+1})>-Cmt$, set $k\leftarrow k+1$ and return
	to Step 2.
	\end{enumerate}
\end{algorithm}

By Proposition \ref{prop:lower bound} and Theorem \ref{thm:exact center}, the algorithm terminates in $O(t)$ iterations, and the resulting analytic center $(D_1,D_2)$ on the central path satisfies $\kappa(D_{1}^{1/2}AD_{2}^{-1/2})-\kappa^\ast \leq 2^{-t}$.

\subsection{Potential Reduction with Approximate Analytic Centers}

Given current $\kappa^0$ and associated analytic center $(D_1(\kappa^0),D_2(\kappa^0))$ and step size $\Delta\kappa$, in order to apply the potential reduction
algorithm, we still need to find the new analytic center $(D_1(\kappa^1),D_2(\kappa^1))$
corresponding to $\kappa^1 = \kappa^{0}-\Delta\kappa$. This can be achieved
using a sequence of Newton iterations starting from $(D_1(\kappa^0),D_2(\kappa^0))$ that will converge to $(D_1(\kappa^1),D_2(\kappa^1))$ quadratically. However, it is not necessary to use exact analytic centers to follow the central path. Instead, we can use ``approximate'' analytic centers $(\overline D_1,\overline D_2)$ that stay close to the central path, and whose associated potentials
\begin{align*}
   \log\det(A^{T}\overline D_{1}A-\overline D_{2})+\log\det(\kappa \overline D_{2}-A^{T}\overline D_{1}A)+\log\det(\overline D_{1}-I_{m})+\log\det(\hat \kappa I_{m}-\overline D_1) 
\end{align*}
are also close to $P(\kappa)$. As we will see, if we start from such approximate
analytic centers, a \emph{single} Newton step with an appropriate step size $\Delta\kappa$
yields a new approximate center that can still achieve a constant reduction in potential function, just like the exact analytic centers in Theorem \ref{thm:exact center}. The advantage is that instead of running Newton updates to convergence at each step, only one Newton update is required. This yields a better computation complexity.
We develop these results in this section. 

To lessen the burden of notation and facilitate understanding of the general approach, we will present our results in this section for the special case of \emph{right-sided} optimal diagonal preconditioning. The generalization to the two-sided problem requires routine adaptations with mostly notational changes. With the optimal right-sided problem, we may set $D_1\equiv I_m$ in the problem $\eqref{eq:SDP-two-sided-upper-bound}$ which then simplifies to 
\begin{align*}
    \min_{D_2 \succeq 0} \kappa \\
M\succeq D_{2}\\
\kappa D_{2}\succeq M
\end{align*}
where $M:=A^TA$. Note that the constraints imply $D_2\succ0,\kappa>0$, so we do not need to explicitly enforce them. With $D$ denoting $D_2$ from now on, the potential function then simplifies to 
\begin{align}
\label{eq:potential-simplified}
P(\kappa):=\max_{D\in \Gamma(\kappa)} & \log\det(M-D)+\log\det(\kappa D-M)+\log\det(D)
\end{align}
Now we focus on developing a potential reduction algorithm based on approximate analytic centers. First, we show that the exact analytic center satisfies an equality condition that is useful in our analysis later. It is essentially the first order condition of the convex optimization problem in \eqref{eq:potential-simplified}. Since it also serves an important role in the Newton step, we state it in a particular form in the next lemma. 
\begin{lemma}
\label{lem:foc}
	The exact analytic center $D$ of the feasible region $\Gamma(\kappa)$
	satisfies 
	\begin{align*}
	Z+\kappa Y & =X
	\end{align*}
	where with $R=M-D$, $S=\kappa D-M$,  \begin{align*}
	X&=R^{-1},Y=S^{-1},Z=D^{-1}.
	\end{align*}
\end{lemma}

By an \textbf{approximate analytic center}, we mean a set of variables $(R,S,D,X,Y,Z)$
that satisfy the optimality conditions 
\begin{align}
\label{eq:optimality}
\begin{split}
R & =M-D\succeq0\\
S & =\kappa D-M\succeq0\\
D & \succeq0\\
Z+\kappa Y & =X,
\end{split}
\end{align}
but instead of requiring $RX=SY=DZ=I$ as is true for the exact analytic
center, we allow the products to be ``close'' to the identity $I$, where we use the following measure of proximity:
\begin{align*}
\delta(A,B) & :=\sqrt{Tr(B^{1/2}AB^{1/2}-I)^{2}}\\
& =\|B^{1/2}AB^{1/2}-I\|_{F}=\|A^{1/2}BA^{1/2}-I\|_{F}
\end{align*}
The second equality above follows from the fact that the Frobenius norm
is a spectral property for symmetric matrices. We may also check that
\begin{align*}
\delta(A,B) & =\sqrt{Tr((AB-I)^{2})}=\sqrt{Tr((BA-I)^{2})}
\end{align*}
As we will see, $\delta(A,B)$ is the appropriate measure
of proximity of $AB$ to $I$. The first result we show is that, starting from the \emph{exact} analytic
center of $\Gamma(\kappa^{0})$, we can obtain a feasible point in
$\Gamma(\kappa^{1})$ that is an $O(\beta)$ \emph{approximate} analytic center of $\Gamma(\kappa^1)$, where $\beta$ is a small constant.
\begin{proposition}
	\label{prop:initial}
	Let $D^{0}:=D(\kappa^0)$ be the exact analytic center of $\Gamma(\kappa^{0})$. Define
	$R^{0}=M-D^{0}$, $S^{0}=\kappa^{0}D^{0}-M$ and $X^{0}=(R^{0})^{-1},Y^{0}=(S^{0})^{-1},Z^{0}=(D^{0})^{-1}$.	Further, let $\kappa^{1}=\kappa^{0}-\Delta\kappa$ with $\Delta\kappa=\frac{\beta}{Tr(D^{0}\cdot(\kappa^{0}D^{0}-M)^{-1})}$
	where $\beta\in(0,1)$ is a sufficiently small constant that does not depend on $D^0,\kappa^0$, and 
	\begin{align*}
	\overline{R}^{0} & =R^{0},\overline{S}^{0}=-\Delta\kappa D^{0}+S^{0},\overline{D}^{0}=D^{0}\\
	\overline{X}^{0} & =X^{0},\overline{Y}^{0}=Y^{0},\overline{Z}^{0}=\Delta\kappa Y^{0}+Z^{0}
	\end{align*}
	Then, the new variables $(\overline{R}^{0},\overline{S}^{0},\overline{D}^{0},\overline{X}^{0},\overline{Y}^{0},\overline{Z}^{0})$ satisfy 
	\begin{align*}
	\overline{R}^{0} & =M-\overline{D}^{0}\succeq0\\
	\overline{S}^{0} & =\kappa^{1}\overline{D}^{0}-M\succeq0\\
	\overline{D}^{0} & \succeq0
	\end{align*}
	and $\overline{Z}^{0}+\kappa^{1}\overline{Y}^{0}=\overline{X}^{0}$.
	Moreover, $\delta(\overline{S}^{0},\overline{Y}^{0})\leq\beta,\delta(\overline{D}^{0},\overline{Z}^{0})\leq\beta,\delta(\overline{R}^{0},\overline{X}^{0})=0$.
\end{proposition}

Next, we show that this set of variables $(\overline{R}^{0},\overline{S}^{0},\overline{D}^{0},\overline{X}^{0},\overline{Y}^{0},\overline{Z}^{0})$
can be used as a starting point to generate a new $ O(\beta^{2}) $ approximate center
for $\Gamma(\kappa^{1})$, using a Newton update. This then implies that starting from an approximate analytic center we can obtain an exact analytic center by applying a sequence of Newton iterations with quadratic convergence. Here to ensure the symmetry
of updates in the Newton step, we use the Nesterov-Todd direction
\citep{todd1998nesterov,nesterov1997self}. 

\begin{proposition}
   \label{prop:Newton step}
   Let $(\overline{R}^{0},\overline{S}^{0},\overline{D}^{0},\overline{X}^{0},\overline{Y}^{0},\overline{Z}^{0})$
	be any set that satisfies the optimality conditions in $\eqref{eq:optimality}$ with $\kappa=\kappa^1$, as well as the proximity conditions
	\begin{align*}
	\delta(\overline{R}^{0},\overline{X}^{0})=\|(\overline{R}^{0})^{\frac{1}{2}}\overline{X}^{0}(\overline{R}^{0})^{\frac{1}{2}}-I\|_{F} & \leq\beta\\
	\delta(\overline{D}^{0},\overline{Z}^{0})=\|(\overline{D}^{0})^{\frac{1}{2}}\overline{Z}^{0}(\overline{D}^{0})^{\frac{1}{2}}-I\|_{F} & \leq\beta\\
	\delta(\overline{S}^{0},\overline{Y}^{0})=\|(\overline{S}^{0})^{\frac{1}{2}}\overline{Y}^{0}(\overline{S}^{0})^{\frac{1}{2}}-I\|_{F} & \leq\beta
	\end{align*}
	for some $\beta\in(0,1)$.	Let $X^{1},Y^{1},Z^{1},R^{1},S^{1},D^{1}$ be defined by 
	\begin{align*}
	X^{1} & =\overline{X}^{0}+\Delta X,Y^{1}=\overline{Y}^{0}+\Delta Y,Z^{1}=\overline{Z}^{0}+\Delta Z\\
	R^{1} & =\overline{R}^{0}+\Delta R,S^{1}=\overline{S}^{0}+\Delta S,D^{1}=\overline{D}^{0}+\Delta D
	\end{align*}
	where the increments $\Delta$s are given by the Newton step 
	\begin{align*}
	\Delta Z+W^{-1}\Delta DW^{-1} & =(\overline{D}^{0})^{-1}-\overline{Z}^{0}\\
	\Delta X+U^{-1}\Delta RU^{-1} & =(\overline{R}^{0})^{-1}-\overline{X}^{0}\\
	\Delta Y+V^{-1}\Delta S  V^{-1}&=(\overline{S}^{0})^{-1}-\overline{Y}^{0}
	\end{align*}
	and 
	$\Delta R  =-\Delta D,\Delta S=\kappa^{1}\Delta D,\Delta X=\Delta Z+\kappa^{1}\Delta Y$,
	where $U,V,W$ are the geometric means 
	\begin{align*}
	U & =(\overline{R}^{0})^{1/2}((\overline{R}^{0})^{1/2}\overline{X}^{0}(\overline{R}^{0})^{1/2})^{-1/2}(\overline{R}^{0})^{1/2}\\
	V & =(\overline{S}^{0})^{1/2}((\overline{S}^{0})^{1/2}\overline{Y}^{0}(\overline{S}^{0})^{1/2})^{-1/2}(\overline{S}^{0})^{1/2}\\
	W & =(\overline{D}^{0})^{1/2}((\overline{D}^{0})^{1/2}\overline{Z}^{0}(\overline{D}^{0})^{1/2})^{-1/2}(\overline{D}^{0})^{1/2}
	\end{align*}
	to ensure that updates are symmetric. Then we have 
	\begin{align*}
	\delta(R^{1},X^{1})=\|(R^{1})^{\frac{1}{2}}X^{1}(R^{1})^{\frac{1}{2}}-I\|_{F} & \leq\frac{1}{2}\frac{\beta^{2}}{1-\beta}\\
	\delta(S^{1},Y^{1})=\|(S^{1})^{\frac{1}{2}}Y^{1}(S^{1})^{\frac{1}{2}}-I\|_{F} & \leq\frac{1}{2}\frac{\beta^{2}}{1-\beta}\\
	\delta(D^{1},Z^{1})=\|(D^{1})^{\frac{1}{2}}Z^{1}(D^{1})^{\frac{1}{2}}-I\|_{F} & \leq\frac{1}{2}\frac{\beta^{2}}{1-\beta}
	\end{align*}
	and 
	\begin{align*}
	\|(\overline{D}^{0})^{-\frac{1}{2}}D^{1}(\overline{D}^{0})^{-\frac{1}{2}}-I\|_{F} & \le\frac{\beta}{1-\beta}.
	\end{align*}
\end{proposition}
One can check that $(R^{1},S^{1},D^{1},X^{1},Y^{1},Z^{1})$ again satisfies the optimality conditions in \eqref{eq:optimality} with $\kappa =\kappa^1$. Note also that the regular Newton step is given by
	\begin{align*}
	\overline{D}^{0}\cdot\Delta Z+\overline{Z}^{0}\cdot\Delta D & =I-\overline{D}^{0}\cdot\overline{Z}^{0}\\
	\overline{R}^{0}\cdot\Delta X+\overline{X}^{0}\cdot\Delta R & =I-\overline{R}^{0}\cdot\overline{X}^{0}\\
	\overline{S}^{0}\cdot\Delta Y+\overline{Y}^{0}\cdot\Delta S & =I-\overline{S}^{0}\cdot\overline{Y}^{0}
	\end{align*}
	but here we need to use the geometric means to ensure the symmetry
	of updates.
	
Our next result shows that starting from an $O(\delta)$ approximate
center $(R^{0},S^{0},D^{0},X^{0},Y^{0},Z^{0})$ instead of the exact analytic
center of $\Gamma(\kappa^{0})$, by following the procedures described in
Proposition \ref{prop:initial} to form the initial set $(\overline{R}^{0},\overline{S}^{0},\overline{D}^{0},\overline{X}^{0},\overline{Y}^{0},\overline{Z}^{0})$
for $\kappa^{1}$, this set of variables satisfies the closeness conditions
\begin{align*}
\delta(\overline{R}^{0},\overline{X}^{0}) & \leq\delta',\delta(\overline{D}^{0},\overline{Z}^{0})\leq\delta',\delta(\overline{S}^{0},\overline{Y}^{0})\leq\delta'
\end{align*}
for $\delta'=\beta+\delta+\beta\delta$, i.e., it is an $O(\delta')$ approximate
analytic center for $\Gamma(\kappa^{1})$.
As long as $\delta$ is small enough, this $O(\delta')$ approximate center can then be used in the
Newton step in Proposition \ref{prop:Newton step} to obtain
an $O((\delta')^{2})=O(\delta)$ approximate center for $\Gamma(\kappa^{1})$,
thus ensuring that we can go from an $O(\delta)$ approximate center
for $\Gamma(\kappa^{0})$ to an $O(\delta)$ approximate center for $\Gamma(\kappa^{1})$
with a single Newton step. 
\begin{proposition}\label{prop:approximate}
	Suppose that $(R^{0},S^{0},D^{0},X^{0},Y^{0},Z^{0})$ satisfy the optimality conditions \eqref{eq:optimality} with $\kappa=\kappa^0$, and that $\delta(R^{0},X^{0})\leq\delta,\delta(D^{0},Z^{0})\leq\delta,\delta(S^{0},Y^{0})\leq\delta$
	for some $\delta\in(0,1)$. Furthermore, let $\kappa^{1}=\kappa^{0}-\Delta\kappa$
	with $
	\Delta\kappa =\frac{\beta}{Tr(D^{0}\cdot(\kappa^{0}D^{0}-M)^{-1})}$ as before, and define
	\begin{align}
	\label{eq:initial-set}
	\begin{split}
	\overline{R}^{0} & =R^{0},\overline{S}^{0}=-\Delta\kappa D^{0}+S^{0},\overline{D}^{0}=D^{0}\\
	\overline{X}^{0} & =X^{0},\overline{Y}^{0}=Y^{0},\overline{Z}^{0}=\Delta\kappa Y^{0}+Z^{0}.
	\end{split}
	\end{align}
	Then we have 
	$\overline{R}^{0} =M-\overline{D}^{0}\succeq0$, $\overline{S}^{0} =\kappa^{1}\overline{D}^{0}-M\succeq 0$, $\overline{D}^{0} \succeq0$, $\overline{Z}^{0}+\kappa^{1}\overline{Y}^{0} =\overline{X}^{0}$,
	and $\delta(\overline{R}^{0},\overline{X}^{0})  \leq\delta$, $\delta(\overline{S}^{0},\overline{Y}^{0})\leq\delta+\beta+\delta\beta$, $\delta(\overline{D}^{0},\overline{Z}^{0})\leq\delta+\beta+\delta\beta$.
\end{proposition}
Propositions \ref{prop:initial}, \ref{prop:Newton step}, and \ref{prop:approximate} guarantee that if we start from an $O(\delta)$
approximate center for $\Gamma(\kappa^{0})$, by forming the initial set according to \eqref{eq:initial-set} with a small step $\Delta\kappa$, and then applying a Newton update step, we can obtain an $O(\delta)$
approximate center for $\kappa^{1}=\kappa^{0}-\Delta\kappa$, and
we may iteratively generate approximate centers in this fashion. In
the next theorem, we show that an approximate center, as the name suggests,
yields a potential value that is $O(\delta^2)$ close to that of the exact center. This will then imply that by constructing a sequence of approximate analytic centers, we can again achieve constant reduction in potential at every step.
\begin{theorem}
	\label{thm:approximate center}Let $\overline{D}^{0}$ and $\overline{D}^{1}$
	be approximate centers for $\kappa^{0}$ and $\kappa^{1}=\kappa^{0}-\frac{\beta}{Tr(\overline{D}^{0}\cdot(\kappa^{0}\overline{D}^{0}-M)^{-1})}$
	respectively, obtained using \eqref{eq:initial-set}
	to generate a feasible point and then applying the Newton update in Proposition \ref{prop:Newton step}. If $	\delta(\overline{R}^{0},\overline{X}^{0}) \leq\delta,\delta(\overline{D}^{0},\overline{Z}^{0})\leq\delta,\delta(\overline{S}^{0},\overline{Y}^{0})\leq\delta$
	for $\delta$ a small universal constant, then 
	\begin{align*}
	\delta(\overline{R}^{1},\overline{X}^{1})  \leq\delta,\delta(\overline{D}^{1},\overline{Z}^{1}) & \leq\delta,\delta(\overline{S}^{1},\overline{Y}^{1})\leq\delta\\
	P(\overline{D}^1,\kappa^1)-P(\kappa^1) &= O(\delta^2)
	\end{align*}
	and $P(\overline{D}^{1},\kappa^{1})\leq P(\overline{D}^{0},\kappa^{0})-c\beta$
	for $c\in(0,1)$ a universal constant, where $P(D,\kappa)$ is defined as 
	\begin{align*}
	   P(D,\kappa):= \log\det(M-D)+\log\det(\kappa D-M)+\log\det(D).
	\end{align*}
\end{theorem}
Theorem \ref{thm:approximate center} allows us to design the following potential algorithm using approximate analytic centers with a similar termination criterion, thanks to the fact that $P(D,\kappa)$ is $O(\delta^2)$ close to $P(\kappa)$.
\begin{algorithm}[Potential Reduction with Approximate Analytic Centers]
	\label{alg:approximate}
	\begin{enumerate}
	\normalfont
	    \item  Let $\overline{D}^{0}$ be an $O(\delta)$ approximate analytic
	center of $\Gamma(\kappa^{0})$ with $\kappa^{0}>\kappa^{\ast}$.
	Set $k=0$. 
	
	\item Let $\Delta\kappa^{k}=\frac{\beta}{Tr(\overline{D}^{0}\cdot(\kappa^{0}\overline{D}^{0}-M)^{-1})}$
	for fixed constant $0<\beta<1$ and let $\kappa^{k+1}=\kappa^{k}-\Delta\kappa^{k}$.
	Update 
	\begin{align*}
	\overline{S}^{0}\leftarrow & -\Delta\kappa\overline{D}^{0}+\overline{S}^{0}\\
	\overline{Z}^{0}\leftarrow & \Delta\kappa Y^{0}+\overline{Z}^{0}
	\end{align*}
	to obtain feasible point for $\kappa^{k+1}$. 
	
	\item Use Newton's update to generate an $O(\delta)$ approximate analytic center
	$\overline{D}^{k+1}$ of $\Gamma(\kappa^{k+1})$. 
	
	\item If $P(\overline{D}^{k+1},\kappa^{k+1})>-Cmt$, set $k\leftarrow k+1$ and return
	to Step 2. 
	\end{enumerate}
\end{algorithm}
Although the potential reduction algorithms presented in this paper are tailored to the optimal diagonal preconditioning problem, we believe that the general approach can be adapted to other quasi-convex problems involving semidefinite matrix variables.
\section{Proofs of Results}
\label{sec:proofs}
\subsection{Proof of Proposition \ref{prop:monotonicity}}
\proof{Proof of Proposition \ref{prop:monotonicity}.}
Let $(D_1(\kappa^0),D_2(\kappa^0))$ and $(D_1(\kappa^1),D_2(\kappa^1))$ denote the analytic centers of $\Gamma(\kappa^{0})$
and $\Gamma(\kappa^{1})$, respectively. We know that $(D_1(\kappa^1),D_2(\kappa^1))$ is
contained in the feasible region $\Gamma(\kappa^{0})$, so that
\begin{align*}
P(\kappa_{0})&\geq\sum_{i}\log\lambda_{i}(A^{T}D_{1}(\kappa^1)A-D_{2}(\kappa^1))+\sum_{i}\log\lambda_{i}(\kappa^0 D_{2}(\kappa^1)-A^{T}D_{1}(\kappa^1)A)\\
&+\sum_{i}\log \lambda_i(D_{1}(\kappa^1)-I_{m})
+\sum_i \log \lambda_i(\hat \kappa I_{m}-D_1(\kappa^1))\\
&>\sum_{i}\log\lambda_{i}(A^{T}D_{1}(\kappa^1)A-D_{2}(\kappa^1))+\sum_{i}\log\lambda_{i}(\kappa^1 D_{2}(\kappa^1)-A^{T}D_{1}(\kappa^1)A)\\
&+\sum_{i}\log \lambda_i(D_{1}(\kappa^1)-I_{m})
+\sum_i \log \lambda_i(\hat \kappa I_{m}-D_1(\kappa^1))\\&=P(\kappa^{1})
\end{align*}
where the strict inequality follows from the Courant-Fischer minimax principle:
for any real symmetric matrix $A$ with $\lambda_{1}\geq\lambda_{2}\dots$, 
\begin{align*}
\lambda_{k} & =\min_{C}\max_{\|x\|=1,Cx=0}\langle Ax,x\rangle
\end{align*}
where $C$ is any $(k-1)\times n$ matrix. This implies that if $A\succeq B$,
$\lambda_{k}(A)\geq\lambda_{k}(B)$, since for 
\begin{align*}
C^\ast & =\arg\min_{C}\max_{\|x\|=1,Cx=0}\langle Ax,x\rangle
\end{align*}
we have 
\begin{align*}
\lambda_{k}(B)\leq\max_{\|x\|=1,C^\ast x=0}\langle Bx,x\rangle & \leq\max_{\|x\|=1,C^\ast x=0}\langle Ax,x\rangle=\lambda_{k}(A)
\end{align*}
and if $A\neq B$, then at least one inequality is strict. 
	\hfill \halmos
\endproof
\subsection{Proof of Proposition \ref{prop:lower bound}}
\proof{Proof of Proposition \ref{prop:lower bound}.}
	Recall that the potential function is given by 
	\begin{align*}
P(\kappa)=\max_{D_{1},D_{2}\in \Gamma(\kappa)} & \log\det(A^{T}D_{1}A-D_{2})+\log\det(\kappa D_{2}-A^{T}D_{1}A)+\log\det(D_{1}-I_{m})+\log\det(\hat \kappa I_{m}-D_1)
\end{align*}
	Since $\overline{\kappa}>\kappa^{\ast}$, $\Gamma(\overline{\kappa})$
	has a bounded and nonempty interior. We claim that there exists an
	interior point $(\overline{D}_1,\overline{D}_2)\in\Gamma(\overline{\kappa})$ such that
	\begin{align*}
	\overline{D}_1-I & \succeq2^{-Ct}I\\
	\hat \kappa I - \overline{D}_1 & \succeq2^{-Ct}I
	\\
	A^T\overline{D}_1A-\overline{D}_2 & \succeq2^{-Ct}I\\
	\overline{\kappa}\overline{D}_2-A^T\overline{D}_1A & \succeq2^{-Ct}I
	\end{align*}
	for some $C$ that only depends on $A$. We may obtain this $(\overline{D}_1,\overline{D}_2)$
	by modifying any $({D}^\ast_1,{D}^\ast_2)$ that satisfies 
	\begin{align*}
	A^T{D}^\ast_1A-D^{\ast}_2 & \succeq0\\
	\kappa^{\ast}D^{\ast}_2-A^T{D}^\ast_1A & \succeq0\\
	D^{\ast}_1-I &
	\succeq0\\
	\hat \kappa I - D^{\ast}_1 &
	\succeq0
	\end{align*}
	as follows. By scaling $({D}^\ast_1,{D}^\ast_2)$ if necessary, we
	can first find some $C$ independent of $({D}^\ast_1,{D}^\ast_2)$ such that $D_1^{\ast}-I\succeq2^{-Ct}I$, $\hat \kappa I-D_1^{\ast}\succeq2^{-Ct}I$, and $D_2^\ast \succeq 2^{-Ct}I$. Then
	for $\delta>0$ sufficiently close to 0, 
	\begin{align*}
	A^TD_1^\ast A-2^{-\delta t}\cdot D_2^{\ast} & \succeq A^TD_1^\ast A-D^{\ast}_2+(D_2^{\ast}-2^{-\delta t}\cdot D_2^{\ast})\\
	& \succeq(1-2^{-\delta t})2^{-Ct}I\succeq2^{-C't}I
	\end{align*}
	Finally, since $\overline{\kappa}>\kappa^{\ast}+2^{-t}$, we have
	\begin{align*}
	\overline{\kappa}2^{-\delta t}D_2^{\ast}-A^TD_1^\ast A & \succeq(\kappa^{\ast}+2^{-t})2^{-\delta t}D_2^{\ast}-A^TD_1^\ast A\\
	& \succeq2^{-(1+\delta)t}D_2^{\ast}\succeq2^{-C't}I
	\end{align*}
	so that $(\overline{D}_1,\overline{D}_2):=(D_1^\ast,2^{-\delta t}D_2^{\ast})$ satisfies the required
	lower bounds. Plugging  $(\overline{D}_1,\overline{D}_2)$ into the definition of $P(\overline{\kappa})$ yields $P(\overline{\kappa})>-Cmt$, where $C$ is a constant that only depends on $A$. 
	\hfill \halmos
\endproof

\subsection{Proof of Theorem \ref{thm:exact center}}
\proof{Proof of Theorem \ref{thm:exact center}.}
We will first show that the potential function $P(\kappa)$ is locally
concave. Note that
\begin{align*}
P((D_{1}D_{2}),\kappa) & :=\log\det(A^{T}D_{1}A-D_{2})+\log\det(\kappa D_{2}-A^{T}D_{1}A)+\log\det(D_{1}-I_{m})+\log\det(\hat{\kappa}I_{m}-D_{1})
\end{align*}
is \emph{jointly} concave in $((D_{1}D_{2}),\kappa)$ on the convex set $\mathcal{D}(\underline{\kappa})\times(\underline{\kappa},\infty)$
for any fixed $\underline{\kappa}>\kappa^{\ast}$, where $\mathcal{D}(\kappa)$
is the interior of the feasible region $\Gamma(\kappa)$. This implies
that, for all $\kappa>\underline{\kappa}$, the following function
\begin{align*}
\underline{P}(\kappa) & :=\max_{(D_{1},D_{2})\in\mathcal{D}(\underline{\kappa})}\log\det(A^{T}D_{1}A-D_{2})+\log\det(\kappa D_{2}-A^{T}D_{1}A)+\log\det(D_{1}-I_{m})+\log\det(\hat{\kappa}I_{m}-D_{1})
\end{align*}
 is a concave function on $(\underline{\kappa},\infty)$. Note, however,
that $\underline{P}(\kappa)$ may not always be equal to the potential
function $P(\kappa)$, because for any $\kappa>\underline{\kappa}$,
the analytic center $(D_{1}(\kappa),D_{2}(\kappa))$ which solves
the problem 
\begin{align*}
\max_{(D_{1},D_{2})\in\mathcal{D}(\kappa)}\log\det(A^{T}D_{1}A-D_{2})+\log\det(\kappa D_{2}-A^{T}D_{1}A)+\log\det(D_{1}-I_{m})+\log\det(\hat{\kappa}I_{m}-D_{1})
\end{align*}
is not necessarily inside $\mathcal{D}(\underline{\kappa})$, as $\underline{\kappa}D_{2}(\kappa)-A^{T}D_{1}(\kappa)A\succ0$
may fail. Fortunately, since the analytic center $(D_{1}(\kappa),D_{2}(\kappa))$
corresponding to each $\kappa>\underline{\kappa}$ is always in the
interior $\mathcal{D}(\kappa)$, with $\kappa^{1}=\kappa^{0}-\Delta\kappa$,
we can guarantee that $(D_{1}(\kappa^{0}),D_{2}(\kappa^{0}))\in\mathcal{D}(\kappa^{1})$
as long as $\Delta\kappa>0$ is small enough. We can then define a
``local'' version of $\underline{P}(\kappa)$ which would coincide
with $P(\kappa)$ and is also concave.

More precisely, at a fixed $\kappa^{0}>\underline{\kappa}$ set $\Delta\kappa=\frac{\beta}{Tr(D_{2}(\kappa^{0})\cdot(\kappa^{0}D_{2}(\kappa^{0})-A^{T}D_{1}(\kappa^{0})A)^{-1})}$
and $\kappa^{1}=\kappa^{0}-\Delta\kappa$, where $\beta$ is a small
constant that is independent of $\kappa^{0}$ to be determined. With
$\lambda_{k}$ denoting the $k$-th largest eigenvalue of a symmetric
matrix, 
\begin{align*}
&|\lambda_{n}(\kappa^{1}D_{2}(\kappa^{0})-A^{T}D_{1}(\kappa^{0})A)-\lambda_{n}(\kappa^{0}D_{2}(\kappa^{0})-A^{T}D_{1}(\kappa^{0})A)| \\
&\leq\Delta\kappa\cdot\max_{i}(D_{2}(\kappa^{0}))_{ii}\\
 & =\frac{\beta}{Tr(D_{2}(\kappa^{0})\cdot(\kappa^{0}D_{2}(\kappa^{0})-A^{T}D_{1}(\kappa^{0})A)^{-1})}\cdot\max_{i}(D_{2}(\kappa^{0}))_{ii}\\
 & \leq\frac{\beta}{\min_{i}(D_{2}(\kappa^{0}))_{ii}\cdot\lambda_{1}((\kappa^{0}D_{2}(\kappa^{0})-A^{T}D_{1}(\kappa^{0})A)^{-1})}\cdot\max_{i}(D_{2}(\kappa^{0}))_{ii}\\
 & =\beta\cdot\kappa(D_{2}(\kappa^{0}))\cdot\lambda_{n}(\kappa^{0}D_{2}(\kappa^{0})-A^{T}D_{1}(\kappa^{0})A)
\end{align*}
 where we have used Weyl's inequality that $\lambda_{n}(X+Y)\leq\lambda_{n}(X)+\lambda_{1}(Y)$
for any symmetric $X,Y$ and 
\begin{align*}
Tr(D_{2}(\kappa^{0})\cdot(\kappa^{0}D_{2}(\kappa^{0})-A^{T}D_{1}(\kappa^{0})A)^{-1}) & \geq Tr(\min_{i}(D_{2}(\kappa^{0}))_{ii}\cdot I\cdot(\kappa^{0}D_{2}(\kappa^{0})-A^{T}D_{1}(\kappa^{0})A)^{-1})\\
 & \geq\min_{i}(D_{2}(\kappa^{0}))_{ii}\cdot\lambda_{1}((\kappa^{0}D_{2}(\kappa^{0})-A^{T}D_{1}(\kappa^{0})A)^{-1})
\end{align*}
For any $\kappa^{0}\in(\kappa^{\ast},\overline{\kappa})$, $\kappa(D_{2}(\kappa^{0}))$
is bounded above by some universal constant that only depends on $A,\overline{\kappa},\hat{\kappa}$.
Therefore, we can make $\beta$ small enough and independent of $\kappa^{0}$
so that $\beta\cdot\kappa(D_{2}(\kappa^{0}))<1$, and 
\begin{align*}
\lambda_{n}(\kappa^{1}D_{2}(\kappa^{0})-A^{T}D_{1}(\kappa^{0})A) & \geq\lambda_{n}(\kappa^{0}D_{2}(\kappa^{0})-A^{T}D_{1}(\kappa^{0})A)-\beta\cdot\kappa(D_{2}(\kappa^{0}))\cdot\lambda_{n}(\kappa^{0}D_{2}(\kappa^{0})-A^{T}D_{1}(\kappa^{0})A)\\
 & >0
\end{align*}
 so that $D_{2}(\kappa^{0})$ is in $\mathcal{D}(\kappa^{1})$. We
have thus proved that, on any interval $[\kappa^{1},\kappa^{0}]\subset(\kappa^{\ast},\overline{\kappa})$
with length bounded by $\frac{\beta}{Tr(D_{2}(\kappa^{0})\cdot(\kappa^{0}D_{2}(\kappa^{0})-A^{T}D_{1}(\kappa^{0})A)^{-1})}$,
$D_{2}(\kappa^{0})$ is in $\mathcal{D}(\kappa^{1})$. Now consider
the function 
\begin{align*}
P_{\kappa^{1}}(\kappa) & :=\max_{(D_{1},D_{2})\in\mathcal{D}(\kappa^{1})}P((D_{1}D_{2}),\kappa)
\end{align*}
 It is concave and coincides with $P(\kappa)$ on $[\kappa^{1},\kappa^{0}]$,
because $D_{2}(\kappa)$ is in $\mathcal{D}(\kappa^{1})$ for all
$\kappa\in[\kappa^{1},\kappa^{0}]$. Therefore, the potential function 
\begin{align*}
P(\kappa)=\max_{(D_{1},D_{2})\in\mathcal{D}(\kappa)}\log\det(A^{T}D_{1}A-D_{2})+\log\det(\kappa D_{2}-A^{T}D_{1}A)+\log\det(D_{1}-I_{m})+\log\det(\hat{\kappa}I_{m}-D_{1})
\end{align*}
 is concave in $\kappa$ on $[\kappa^{1},\kappa^{0}]$. 

Next we show that the potential function on the interval $[\kappa^{1},\kappa^{0}]$
increases by a constant lower bound. Note that
\begin{align*}
P(\kappa) & =P((D_{1}(\kappa),D_{2}(\kappa)),\kappa)
\end{align*}
 where the analytic center $(D_{1}(\kappa),D_{2}(\kappa))$ is the
unique maximizer corresponding to $\kappa$, so that 
\begin{align*}
\frac{d}{d\kappa}P(\kappa) & =Tr(\nabla_{(D_{1},D_{2})}P((D_{1}(\kappa),D_{2}(\kappa)),\kappa)\cdot\partial_{\kappa}(D_{1}(\kappa),D_{2}(\kappa)))+\partial_{\kappa}P((D_{1}(\kappa),D_{2}(\kappa)),\kappa)\\
 & =\partial_{\kappa}P((D_{1}(\kappa),D_{2}(\kappa)),\kappa)\\
 & =Tr(D_{2}(\kappa)\cdot(\kappa D_{2}(\kappa)-A^{T}D_{1}(\kappa)A)^{-1})
\end{align*}
 using the envelope theorem, or equivalently the first order condition
of the maximization problem in the definition of the potential function. By concavity of $P(\kappa)$, we have
\begin{align*}
P(\kappa^{1})-P(\kappa^{0}) & \leq\frac{d}{d\kappa}P(\kappa)\mid_{\kappa^{0}}\cdot(\kappa^{1}-\kappa^{0})\\
 & =-\Delta\kappa\cdot Tr(D_{2}(\kappa^{0})\cdot(\kappa^{0}D_{2}(\kappa^{0})-A^{T}D_{1}(\kappa^{0})A)^{-1})\\
 & =-\beta
\end{align*}
and since $P(\kappa)$ is monotone increasing, we have proved that
for any $\kappa^{0}\in(\kappa^{\ast},\overline{\kappa})$, decreasing
$\kappa^{0}$ by 
\begin{align*}
\Delta\kappa & =\frac{\beta}{Tr(D_{2}(\kappa^{0})\cdot(\kappa^{0}D_{2}(\kappa^{0})-A^{T}D_{1}(\kappa^{0})A)^{-1})}
\end{align*}
 is guranteed to \emph{decrease} the potential function by at least
$\beta$. 
	\hfill \halmos
\endproof

\subsection{Proof of Lemma \ref{lem:foc}}
\proof{Proof of Lemma \ref{lem:foc}.}
	Since $D$ solves the convex problem
	\begin{align*}
	& \max_{D}\log\det(M-D)+\log\det(\kappa D-M)+\log\det D
	\end{align*}
	The first order condition implies 
	\begin{align*}
	-(M-D)^{-1}+\kappa(\kappa D-M)^{-1}+D^{-1} & =0
	\end{align*}
	so if we define $R=M-D$ and $S=\kappa D-M$, and 
	\begin{align*}
	X	=R^{-1},Y=S^{-1},Z=D^{-1}
	\end{align*}
	then we have $	Z+\kappa Y=X$ and moreover $R+S=(\kappa-1)D$. 
	\hfill \halmos
\endproof

\subsection{Proof of Proposition \ref{prop:initial}}
\proof{Proof of Proposition \ref{prop:initial}.}
	We can readily verify that the new variables satisfy the desired equalities.
	When $\Delta\kappa$ is appropriately chosen as in Theorem \ref{thm:exact center},
	we can guarantee that $\overline{S}^{0}\succeq0$. Now 
	\begin{align*}
	\delta(\overline{S}^{0},\overline{Y}^{0})^{2}=Tr((\overline{S}^{0}\overline{Y}^{0}-I)^{2}) & =Tr((-\Delta\kappa D^{0}+S^{0})Y^{0}-I)^{2}\\
	& =Tr(-\Delta\kappa D^{0}Y^{0})^{2}\\
	& =Tr(-\Delta\kappa D^{0}(\kappa^{0}D^{0}-M)^{-1})^{2}\\
	& \leq(Tr(-\Delta\kappa D^{0}(\kappa^{0}D^{0}-M)^{-1}))^{2}=\beta^{2}
	\end{align*}
	where we have used $Tr(A^{2})=(TrA)^{2}-2\sigma_{2}(A)$ where $\sigma_{2}(A)=\sum_{i<j}\lambda_{i}\lambda_{j}$
	and $\lambda_{i}$ are the eigenvalues of $A$. Therefore, $\delta(\overline{S}^{0},\overline{Y}^{0})\leq\beta$.
	Similarly, $\delta(\overline{D}^{0},\overline{Z}^{0})\leq\beta$ and
	since $\overline{R}^{0}=R^{0},\overline{X}^{0}=X^{0}$, we have $\delta(\overline{D}^{0},\overline{Z}^{0})=0$
	\hfill \halmos
\endproof
\subsection{Proof of Proposition \ref{prop:Newton step}}
\proof{Proof of Proposition \ref{prop:Newton step}.}
	First observe that the Newton system has a unique solution. To see
	this, note that we can eliminate $\Delta X,\Delta Y,\Delta Z$ using
	$\Delta Z+\kappa^{1}\Delta Y=\Delta X$. Then using
	\begin{align*}
	\Delta R & =-\Delta D\\
	\Delta S & =\kappa^{1}\Delta D
	\end{align*}
	we can uniquely solve for $\Delta D$, from which all other increments
	can be uniquely determined. 
	
	Next, we have
	\begin{align*}
	\delta(D^{1},Z^{1})^{2}=\|(Z^{1})^{\frac{1}{2}}D^{1}(Z^{1})^{\frac{1}{2}}-I\|_{F}^{2} & =Tr(((Z^{1})^{\frac{1}{2}}D^{1}(Z^{1})^{\frac{1}{2}}-I)^{2})\\
	& =Tr((Z^{1})D^{1}(Z^{1})D^{1})-2Tr(Z^{1}D^{1})+Tr(I)\\
	& =Tr((D^{1}Z^{1}-I)^{2})\\
	& \leq\left(Tr(D^{1}Z^{1}-I)\right)^{2}
	\end{align*}
	where we have used $Tr(A^{2})=(TrA)^{2}-2\sigma_{2}(A)$ where 
	\begin{align*}
	\sigma_{2}(A) & =\sum_{i<j}\lambda_{i}\lambda_{j}
	\end{align*}
	and $\lambda_{i}$ are the eigenvalues of $A$. 
	
	The Newton update formulae 
	\begin{align*}
	\Delta Z+W^{-1}\Delta DW^{-1} & =(\overline{D}^{0})^{-1}-\overline{Z}^{0}\\
	Z^{1} & =Z^{0}+\Delta Z\\
	D^{1} & =D^{0}+\Delta D
	\end{align*}
	and $W =(\overline{D}^{0})^{1/2}((\overline{D}^{0})^{1/2}\overline{Z}^{0}(\overline{D}^{0})^{1/2})^{-1/2}(\overline{D}^{0})^{1/2}$
	imply
	\begin{align*}
	Tr(D^{1}Z^{1}-I) & =Tr((\overline{D}^{0}+\Delta D)(\overline{Z}^{0}+\Delta Z)-I)\\
	& =Tr(\overline{D}^{0}\overline{Z}^{0}+\Delta D\overline{Z}^{0}+\overline{D}^{0}\Delta Z+\Delta D\Delta Z-I)\\
	& =Tr(-\overline{D}^{0}W^{-1}\Delta DW^{-1}+\Delta D\overline{Z}^{0}+\Delta D\Delta Z)\\
	& =Tr(-(\overline{D}^{0})^{1/2}\overline{Z}^{0}(\overline{D}^{0})^{1/2}(\overline{D}^{0})^{-1/2}\Delta D(\overline{D}^{0})^{-1/2}+\Delta D\overline{Z}^{0}+\Delta D\Delta Z)\\
	& =Tr(-\Delta D\overline{Z}^{0}+\Delta D\overline{Z}^{0}+\Delta D\Delta Z)\\
	& =Tr(\Delta D\Delta Z)
	\end{align*}
	So that 
	\begin{align*}
	\delta(D^{1},Z^{1}) & \leq Tr(\Delta D\Delta Z)
	\end{align*}
	and similarly
	\begin{align*}
	\delta(R^{1},X^{1}) & \leq Tr(\Delta R\Delta X)\\
	\delta(S^{1},Y^{1}) & \leq Tr(\Delta S\Delta Y).
	\end{align*}
	
	Multiplying both sides of
	\begin{align*}
	\Delta Z+W^{-1}\Delta DW^{-1} & =(\overline{D}^{0})^{-1}-\overline{Z}^{0}
	\end{align*}
	by $W^{\frac{1}{2}}$ on the left and on the right we have
	\begin{align*}
	W^{\frac{1}{2}}\Delta ZW^{\frac{1}{2}}+W^{-\frac{1}{2}}\Delta DW^{-\frac{1}{2}} & =W^{\frac{1}{2}}(\overline{D}^{0})^{-1}W^{\frac{1}{2}}-W^{\frac{1}{2}}\overline{Z}^{0}W^{\frac{1}{2}}
	\end{align*}
	Then squaring and taking the trace, 
	\begin{align*}
	Tr(W^{\frac{1}{2}}\Delta ZW^{\frac{1}{2}}+W^{-\frac{1}{2}}\Delta DW^{-\frac{1}{2}})^{2} & =Tr(W^{\frac{1}{2}}(\overline{D}^{0})^{-1}W^{\frac{1}{2}}-W^{\frac{1}{2}}\overline{Z}^{0}W^{\frac{1}{2}})^{2}
	\end{align*}
	and the left hand side is equal to 
	\begin{align*}
	Tr(W^{\frac{1}{2}}\Delta ZW^{\frac{1}{2}}+W^{-\frac{1}{2}}\Delta DW^{-\frac{1}{2}})^{2} & =Tr(W^{\frac{1}{2}}\Delta ZW^{\frac{1}{2}})^{2}+Tr(W^{-\frac{1}{2}}\Delta DW^{-\frac{1}{2}})^{2}\\
	& +Tr(\Delta Z\Delta D)+Tr(\Delta D\Delta Z)
	\end{align*}
	So we have 
	\begin{align*}
	2Tr(\Delta Z\Delta D) & \leq Tr(W^{\frac{1}{2}}(\overline{D}^{0})^{-1}W^{\frac{1}{2}}-W^{\frac{1}{2}}\overline{Z}^{0}W^{\frac{1}{2}})^{2}
	\end{align*}
	Now we calculate the right hand side. Using
	\begin{align*}
	W & =(\overline{D}^{0})^{1/2}((\overline{D}^{0})^{1/2}\overline{Z}^{0}(\overline{D}^{0})^{1/2})^{-1/2}(\overline{D}^{0})^{1/2}\\
	& =(\overline{Z}^{0})^{-1/2}((\overline{Z}^{0})^{1/2}\overline{D}^{0}(\overline{Z}^{0})^{1/2})^{1/2}(\overline{Z}^{0})^{-1/2}
	\end{align*}
	we get 
	\begin{align*}
	Tr(W^{\frac{1}{2}}(\overline{D}^{0})^{-1}W^{\frac{1}{2}})^{2} & =Tr(\overline{D}^{0})^{1/2}((\overline{D}^{0})^{1/2}\overline{Z}^{0}(\overline{D}^{0})^{1/2})^{-1}(\overline{D}^{0})^{1/2}(\overline{D}^{0})^{-1}\\
	& =Tr(\overline{Z}^{0})^{-1}(\overline{D}^{0})^{-1}\\
	Tr(W^{\frac{1}{2}}\overline{Z}^{0}W^{\frac{1}{2}})^{2} & =Tr(\overline{Z}^{0})^{-1/2}((\overline{Z}^{0})^{1/2}\overline{D}^{0}(\overline{Z}^{0})^{1/2})^{1}(\overline{Z}^{0})^{-1/2}\overline{Z}^{0}\\
	& =Tr\overline{D}^{0}\overline{Z}^{0}\\
	Tr(W^{\frac{1}{2}}(\overline{D}^{0})^{-1}W^{\frac{1}{2}})(W^{\frac{1}{2}}\overline{Z}^{0}W^{\frac{1}{2}}) & =Tr(I)\\
	Tr(W^{\frac{1}{2}}\overline{Z}^{0}W^{\frac{1}{2}})(W^{\frac{1}{2}}(\overline{D}^{0})^{-1}W^{\frac{1}{2}}) & =Tr(I)
	\end{align*}
	so that %
	\begin{align*}
	Tr(W^{\frac{1}{2}}(\overline{D}^{0})^{-1}W^{\frac{1}{2}}-W^{\frac{1}{2}}\overline{Z}^{0}W^{\frac{1}{2}})^{2} & =Tr(\overline{Z}^{0})^{-1}(\overline{D}^{0})^{-1}+Tr\overline{D}^{0}\overline{Z}^{0}+2Tr(I)\\
	& =Tr((\overline{D}^{0})^{-1/2}(\overline{Z}^{0})^{-1}(\overline{D}^{0})^{-1/2}(I-2(\overline{D}^{0})^{1/2}\overline{Z}^{0}(\overline{D}^{0})^{1/2}+((\overline{D}^{0})^{1/2}\overline{Z}^{0}(\overline{D}^{0})^{1/2})^{2}))\\
	& =Tr((\overline{D}^{0})^{-1/2}(\overline{Z}^{0})^{-1}(\overline{D}^{0})^{-1/2}(I-(\overline{D}^{0})^{1/2}\overline{Z}^{0}(\overline{D}^{0})^{1/2})^{2})
	\end{align*}
	
	Now recall that 
	\begin{align*}
	\delta(\overline{D}^{0},\overline{Z}^{0})=\|(\overline{D}^{0})^{\frac{1}{2}}\overline{Z}^{0}(\overline{D}^{0})^{\frac{1}{2}}-I\|_{F} & \leq\beta
	\end{align*}
	and so 
	\begin{align*}
	Tr((I-(\overline{D}^{0})^{1/2}\overline{Z}^{0}(\overline{D}^{0})^{1/2})^{2}) & =\|(\overline{D}^{0})^{1/2}\overline{Z}^{0}(\overline{D}^{0})^{1/2}-I\|_{F}^{2}\\
	& =\delta(\overline{D}^{0},\overline{Z}^{0})^{2}\leq\beta^{2}
	\end{align*}
	and H{\"o}lder's inequality yields 
	\begin{align*}
	Tr((\overline{D}^{0})^{-1/2}(\overline{Z}^{0})^{-1}(\overline{D}^{0})^{-1/2}(I-(\overline{D}^{0})^{1/2}\overline{Z}^{0}(\overline{D}^{0})^{1/2})^{2}) & \leq\lambda_{1}((\overline{D}^{0})^{-1/2}(\overline{Z}^{0})^{-1}(\overline{D}^{0})^{-1/2})\cdot Tr(I-(\overline{D}^{0})^{1/2}\overline{Z}^{0}(\overline{D}^{0})^{1/2})^{2})\\
	& \leq\lambda_{1}((\overline{Z}^{0})^{-1}(\overline{D}^{0})^{-1})\beta^{2}
	\end{align*}
	Recall the eigenvalue stability result 
	\begin{align*}
	|\lambda_{i}(X)-\lambda_{i}(Y)| & \leq\|X-Y\|_{op}
	\end{align*}
	for symmetric positive definite matrices $X,Y$ and again using 
	\begin{align*}
	\|(\overline{D}^{0})^{1/2}\overline{Z}^{0}(\overline{D}^{0})^{1/2}-I\|_{F} & =\|\overline{D}^{0}\overline{Z}^{0}-I\|_{F}\\
	& =\sqrt{Tr(\overline{D}^{0}\overline{Z}^{0}-I)^{T}(\overline{D}^{0}\overline{Z}^{0}-I)}\leq\beta
	\end{align*}
	and the relation $\|X\|_{op} \leq\|X\|_{F}$
	between operator norm and Frobenius norm, we have 
	\begin{align*}
	|\lambda_{n}((\overline{D}^{0})^{1/2}\overline{Z}^{0}(\overline{D}^{0})^{1/2})-1| & =|\lambda_{n}((\overline{D}^{0})^{1/2}\overline{Z}^{0}(\overline{D}^{0})^{1/2})-\lambda_{n}(I)|\\
	& \leq\|(\overline{D}^{0})^{1/2}\overline{Z}^{0}(\overline{D}^{0})^{1/2}-I\|_{op}\leq\|(\overline{D}^{0})^{1/2}\overline{Z}^{0}(\overline{D}^{0})^{1/2}-I\|_{F}\leq\beta
	\end{align*}
	so that 
	\begin{align*}
	\lambda_{n}((\overline{D}^{0})^{1/2}\overline{Z}^{0}(\overline{D}^{0})^{1/2}) & \geq1-\beta
	\end{align*}
	and 
	\begin{align*}
	\lambda_{1}((\overline{D}^{0})^{-1/2}(\overline{Z}^{0})^{-1}(\overline{D}^{0})^{-1/2}) & =\frac{1}{\lambda_{n}((\overline{D}^{0})^{1/2}\overline{Z}^{0}(\overline{D}^{0})^{1/2})}\leq\frac{1}{1-\beta}
	\end{align*}
	This allows us to finally conclude that 
	\begin{align*}
	Tr(\Delta Z\Delta D) & \leq\frac{1}{2}\frac{\beta^{2}}{1-\beta}
	\end{align*}
	and so 
	\begin{align*}
	\delta(D^{1},Z^{1}) & \leq Tr(\Delta Z\Delta D)\\
	& \leq\frac{1}{2}\frac{\beta^{2}}{1-\beta}
	\end{align*}
	The same argument applied to $(R^{1},X^{1})$, $(S^{1},Y^{1})$ gives
	\begin{align*}
	\delta(R^{1},X^{1}) & \leq\frac{1}{2}\frac{\beta^{2}}{1-\beta}\\
	\delta(S^{1},Y^{1}) & \leq\frac{1}{2}\frac{\beta^{2}}{1-\beta}
	\end{align*}
	
	Now we show that 
	\begin{align*}
	\|(\overline{D}^{0})^{-\frac{1}{2}}D^{1}(\overline{D}^{0})^{-\frac{1}{2}}-I\|_{F} & \le\frac{\beta}{1-\beta}
	\end{align*}
	Note 
	\begin{align*}
	\|(\overline{D}^{0})^{-\frac{1}{2}}D^{1}(\overline{D}^{0})^{-\frac{1}{2}}-I\|_{F}^{2} & =Tr(((\overline{D}^{0})^{-\frac{1}{2}}D^{1}(\overline{D}^{0})^{-\frac{1}{2}}-I)^{2})\\
	& =Tr(((\overline{D}^{0})^{-1}(D^{1}-\overline{D}^{0})(\overline{D}^{0})^{-1}(D^{1}-\overline{D}^{0})\\
	& =Tr((\overline{D}^{0})^{-1}\Delta D(\overline{D}^{0})^{-1}\Delta D)
	\end{align*}
	Now multiplying both sides of 
	\begin{align*}
	W\Delta ZW+\Delta D=(\overline{Z}^{0})^{-1}-\overline{D}^{0}
	\end{align*}
	by $(\overline{D}^{0})^{-\frac{1}{2}}$, we have 
	\begin{align*}
	Tr(((\overline{D}^{0})^{-1}\Delta D)^{2}) & \leq Tr(((\overline{D}^{0})^{-1}(\overline{Z}^{0})^{-1}-I))^{2})\\
	& =\|(\overline{D}^{0})^{-\frac{1}{2}}(\overline{Z}^{0})^{-1}(\overline{D}^{0})^{-\frac{1}{2}}-I\|_{F}^{2}\\
	& =\|(I-(\overline{D}^{0})^{\frac{1}{2}}(\overline{Z}^{0})(\overline{D}^{0})^{\frac{1}{2}})((\overline{D}^{0})^{-\frac{1}{2}}(\overline{Z}^{0})^{-1}(\overline{D}^{0})^{-\frac{1}{2}})\|_{F}^{2}\\
	& \leq Tr((I-(\overline{D}^{0})^{\frac{1}{2}}(\overline{Z}^{0})(\overline{D}^{0})^{\frac{1}{2}})^{2})\lambda_{1}(((\overline{D}^{0})^{-\frac{1}{2}}(\overline{Z}^{0})^{-1}(\overline{D}^{0})^{-\frac{1}{2}})^{2})\\
	& \leq\beta^{2}/(1-\beta)^{2}
	\end{align*}
	as desired. 
	\hfill \halmos
\endproof
\subsection{Proof of Proposition \ref{prop:approximate}}
\proof{Proof of Proposition \ref{prop:approximate}.}
We have $\delta(\overline{R}^{0},\overline{X}^{0})=\delta(R^{0},X^{0})\leq\delta$
and 
\begin{align*}
\delta(\overline{S}^{0},\overline{Y}^{0})^{2} & =Tr(((\overline{Y}^{0})^{\frac{1}{2}}\overline{S}^{0}(\overline{Y}^{0})^{\frac{1}{2}}-I)^{2})\\
 & =Tr(((Y^{0})^{\frac{1}{2}}(-\Delta\kappa D^{0}+S^{0})(Y^{0})^{\frac{1}{2}}-I)^{2})\\
 & =Tr(((Y^{0})^{\frac{1}{2}}S^{0}(Y^{0})^{\frac{1}{2}}-\Delta\kappa(Y^{0})^{\frac{1}{2}}D^{0}(Y^{0})^{\frac{1}{2}}-I)^{2})\\
 & =\|(Y^{0})^{\frac{1}{2}}S^{0}(Y^{0})^{\frac{1}{2}}-I-\Delta\kappa(Y^{0})^{\frac{1}{2}}D^{0}(Y^{0})^{\frac{1}{2}}\|_{F}^{2}
\end{align*}
 so using the triangle inequality of the Frobenius norm, 
\begin{align*}
\delta(\overline{S}^{0},\overline{Y}^{0}) & =\|(Y^{0})^{\frac{1}{2}}S^{0}(Y^{0})^{\frac{1}{2}}-I-\Delta\kappa(Y^{0})^{\frac{1}{2}}D^{0}(Y^{0})^{\frac{1}{2}}\|_{F}\\
 & \leq\|(Y^{0})^{\frac{1}{2}}S^{0}(Y^{0})^{\frac{1}{2}}-I\|_{F}+\|\Delta\kappa(Y^{0})^{\frac{1}{2}}D^{0}(Y^{0})^{\frac{1}{2}}\|_{F}\\
 & \leq\delta+Tr(\Delta\kappa(Y^{0})^{\frac{1}{2}}D^{0}(Y^{0})^{\frac{1}{2}})
\end{align*}
where the second term can be further bounded by
\begin{align*}
Tr(\Delta\kappa(Y^{0})^{\frac{1}{2}}D^{0}(Y^{0})^{\frac{1}{2}}) & =\Delta\kappa Tr(S^{0}Y^{0}D^{0}(S^{0})^{-1})\\
 & =\Delta\kappa Tr(D^{0}(S^{0})^{-1}))+\Delta\kappa Tr((S^{0}Y^{0}-I)D^{0}(S^{0})^{-1})\\
 & =\beta+\Delta\kappa Tr((S^{0}Y^{0}-I)D^{0}(S^{0})^{-1})\\
 & \leq\beta+\Delta\kappa\sqrt{Tr((S^{0}Y^{0}-I)^{2})}\cdot\sqrt{Tr((D^{0}(S^{0})^{-1})^{2})}\\
 & \leq\beta+\Delta\kappa\cdot\delta\cdot Tr(D^{0}(S^{0})^{-1})=\beta+\delta\beta
\end{align*}
so that $\delta(\overline{S}^{0},\overline{Y}^{0})\leq\delta+\beta+\delta\beta$,
and similarly $\delta(\overline{D}^{0},\overline{Z}^{0})\leq\delta+\beta+\delta\beta$.
\hfill \halmos
\endproof

\subsection{Proof of Theorem \ref{thm:approximate center}}
\proof{Proof of Theorem \ref{thm:approximate center}.}
First we verify that with $\delta$ small enough, the approximate analytic centers for $\kappa^1$ obtained by applying the updates \eqref{eq:initial-set} followed by a Newton update
satisfy $	\delta(\overline{R}^{1},\overline{X}^{1}),\delta(\overline{D}^{1},\overline{Z}^{1}),\delta(\overline{S}^{1},\overline{Y}^{1})=O(\delta)$.
Proposition \ref{prop:approximate} combined with Proposition \ref{prop:Newton step} imply that
	\begin{align*}
	\delta(\overline{R}^{1},\overline{X}^{1})\leq\frac{1}{2}\frac{(\delta+\beta+\delta\beta)^{2}}{1-(\delta+\beta+\delta\beta)} & \leq\delta,
	\end{align*}
as long as $\delta,\beta$ are sufficiently small. 
	Now recall the definition 
	\begin{align*}
	P(D,\kappa) & :=\log\det(M-D)+\log\det(\kappa D-M)+\log\det D
	\end{align*}
	The main part of the proof consists of showing that when an approximate analytic
	center $D$ for $\Gamma(\kappa)$ is $O(\eta)$ close to the exact analytic center, the potential function $P(D,\kappa)$ is
	$O(\eta^{2})$ close to the potential function $P(\kappa)$, for sufficiently small
	$\eta$.

	For fixed $\kappa>\kappa^\ast$, suppose we have an approximate center $\overline{D}$
	with 
	\begin{align*}
	\overline{R} & =M-\overline{D}\succeq0\\
	\overline{S} & =\kappa\overline{D}-M\succeq0\\
	\overline{D} & \succeq0\\
	\overline{Z}+\kappa\overline{Y} & =\overline{X}
	\end{align*}
	and $\delta(\overline{R},\overline{X})  \leq\eta,
	\delta(\overline{S},\overline{Y})  \leq\eta,
	\delta(\overline{D},\overline{Z})  \leq\eta$ for some small $\eta\in(0,1)$. Let $D$ be the exact analytic center of
	$\Gamma(\kappa)$. We will first show that $	\|(\overline{D})^{-\frac{1}{2}}D(\overline{D})^{-\frac{1}{2}}-I\|_{F}=O(\eta)$. Proposition \ref{prop:Newton step} applied to $\overline{D}$
	implies that 
	\begin{align*}
	\|(\overline{D})^{-\frac{1}{2}}D'(\overline{D})^{-\frac{1}{2}}-I\|_{F} & \le\frac{\eta}{1-\eta}
	\end{align*}
	where $D'$ is the result of one step of Newton update applied to
	$\overline{D}$, and is a new $\frac{1}{2}\frac{\eta^{2}}{1-\eta}$ approximate center. 
	Now
	the key is to iterate the Newton update so that we obtain a sequence
	of approximate centers $D^{(1)},D^{(2)},\dots$, that converges to
	$D$. Using the upper bounds on the approximate centers, we obtain,
	with $\overline{D}=D^{(0)}$, 
	\begin{align*}
	\|(\overline{D})^{-\frac{1}{2}}D(\overline{D})^{-\frac{1}{2}}-I\|_{F} & =\|(\overline{D})^{-\frac{1}{2}}D(\overline{D})^{-\frac{1}{2}}-(\overline{D})^{-\frac{1}{2}}(D^{(1)})(\overline{D})^{-\frac{1}{2}}+(\overline{D})^{-\frac{1}{2}}(D^{(1)})(\overline{D})^{-\frac{1}{2}}-I\|_{F}\\
	& =\|(D^{(0)})^{-\frac{1}{2}}D(D^{(0)})^{-\frac{1}{2}}-(D^{(0)})^{-\frac{1}{2}}(D^{(1)})(D^{(0)})^{-\frac{1}{2}}+(D^{(0)})^{-\frac{1}{2}}(D^{(1)})(D^{(0)})^{-\frac{1}{2}}-I\|_{F}\\
	& \leq\|(D^{(0)})^{-\frac{1}{2}}D(D^{(0)})^{-\frac{1}{2}}-(D^{(0)})^{-\frac{1}{2}}(D^{(1)})(D^{(0)})^{-\frac{1}{2}}\|_{F}+\|(D^{(0)})^{-\frac{1}{2}}(D^{(1)})(D^{(0)})^{-\frac{1}{2}}-I\|_{F}\\
	& \leq\sqrt{Tr((D^{(0)})^{-1}(D-(D^{(1)}))^{2}}+\frac{\eta}{1-\eta}\\
	& =\sqrt{Tr((D^{(0)})^{-1}(D^{(1)})(D^{(1)})^{-1}(D-(D^{(1)}))^{2}}+\frac{\eta}{1-\eta}\\
	& \leq\lambda_{1}((D^{(0)})^{-\frac{1}{2}}(D^{(1)})(D^{(0)})^{-\frac{1}{2}})\cdot\sqrt{Tr((D^{(1)})^{-\frac{1}{2}}D(D^{(1)})^{-\frac{1}{2}}-I)^{2}}+\frac{\eta}{1-\eta}\\
	& =\lambda_{1}((D^{(0)})^{-\frac{1}{2}}(D^{(1)})(D^{(0)})^{-\frac{1}{2}})\cdot\|(D^{(1)})^{-\frac{1}{2}}D(D^{(1)})^{-\frac{1}{2}}-I\|_{F}+\frac{\eta}{1-\eta}
	\end{align*}
	where we have used the property $\|AB\|_F \leq \lambda_1(A)\cdot\|B\|_F$
	to obtain the last inequality. Now since 
	\begin{align*}
	\|(D^{(0)})^{-\frac{1}{2}}(D^{(1)})(D^{(0)})^{-\frac{1}{2}}-I\|_{F} & \leq\frac{\eta}{1-\eta}
	\end{align*}
	we have 
	\begin{align*}
	|\lambda_{1}((D^{(0)})^{-\frac{1}{2}}(D^{(1)})(D^{(0)})^{-\frac{1}{2}})| & \leq\lambda_{1}((D^{(0)})^{-\frac{1}{2}}(D^{(1)})(D^{(0)})^{-\frac{1}{2}}-I)+\lambda_{1}(I)\\
	& \leq\|(D^{(0)})^{-\frac{1}{2}}(D^{(1)})(D^{(0)})^{-\frac{1}{2}}-I\|_{F}+1\\
	& \leq\frac{1}{1-\eta}.
	\end{align*}
	Now to bound $\|(D^{(1)})^{-\frac{1}{2}}D(D^{(1)})^{-\frac{1}{2}}-I\|_{F}$, we can use the same argument as above, with the only difference that, since $ D^{(1)}$ is a $\frac{1}{2}\frac{\eta^{2}}{1-\eta}$ approximate center, Proposition \ref{prop:Newton step} applied to $ D^{(1)}$ implies 
	\begin{align*}
	\|(D^{(1)})^{-\frac{1}{2}}D^{(2)}(D^{(1)})^{-\frac{1}{2}}-I\|_{F} & \le\frac{\frac{1}{2}\frac{\eta^{2}}{1-\eta}}{1-\frac{1}{2}\frac{\eta^{2}}{1-\eta}}\leq\eta^2
	\end{align*}
	for $\eta<1/2$, since $\frac{1}{2}\frac{\eta^{2}}{1-\eta}\leq\eta^{2}<\eta$. Applying the argument above to $\|(D^{(1)})^{-\frac{1}{2}}D(D^{(1)})^{-\frac{1}{2}}-I\|_{F}$ then yields
	\begin{align*}
	\|(\overline{D})^{-\frac{1}{2}}D(\overline{D})^{-\frac{1}{2}}-I\|_{F} & \leq\frac{1}{1-\eta}\cdot(\frac{1}{1-\frac{1}{2}\frac{\eta^{2}}{1-\eta}}\cdot\|(D^{(2)})^{-\frac{1}{2}}D(D^{(2)})^{-\frac{1}{2}}-I\|_{F}+\frac{\frac{1}{2}\frac{\eta^{2}}{1-\eta}}{1-\frac{1}{2}\frac{\eta^{2}}{1-\eta}})+\frac{\eta}{1-\eta}\\
	& \leq\frac{1}{1-\eta}\cdot\frac{1}{1-\frac{1}{2}\frac{\eta^{2}}{1-\eta}}\cdot\|(D^{(2)})^{-\frac{1}{2}}D(D^{(2)})^{-\frac{1}{2}}-I\|_{F}+\frac{\eta^{2}}{(1-\eta)^{2}}+\frac{\eta}{1-\eta}\\
	& \leq\frac{1}{1-\eta}\cdot\frac{1}{1-\eta^{2}}\cdot\|(D^{(2)})^{-\frac{1}{2}}D(D^{(2)})^{-\frac{1}{2}}-I\|_{F}+\frac{\eta^{2}}{(1-\eta)^{2}}+\frac{\eta}{1-\eta}
	\end{align*}
	Now we can use the fact that $D^{(k)}$ is a $\eta^{2^{k}}$ approximate
center to prove the induction argument
\begin{align*}
\prod_{i=1}^{k}(\frac{1}{1-\eta{}^{2^{i-1}}})\|(D^{(k)})^{-\frac{1}{2}}D(D^{(k)})^{-\frac{1}{2}}-I\|_{F}+\sum_{i=1}^{k}\frac{\eta^{i}}{(1-\eta)^{i}} & \leq\\
\prod_{i=1}^{k}(\frac{1}{1-\eta{}^{2^{i-1}}})\cdot(\frac{1}{1-\eta^{2^{k}}}\cdot\|(D^{(k+1)})^{-\frac{1}{2}}D(D^{(k+1)})^{-\frac{1}{2}}-I\|_{F}+\frac{\eta^{2^{k}}}{1-\eta^{2^{k}}})+\sum_{i=1}^{k}\frac{\eta^{i}}{(1-\eta)^{i}} & \leq\\
\prod_{i=1}^{k+1}(\frac{1}{1-\eta{}^{2^{i-1}}})\cdot\|(D^{(k+1)})^{-\frac{1}{2}}D(D^{(k+1)})^{-\frac{1}{2}}-I\|_{F}+\sum_{i=1}^{k+1}\frac{\eta^{i}}{(1-\eta)^{i}}
\end{align*}
 and since $\lim_{k\rightarrow\infty}\prod_{i=1}^{k}(\frac{1}{1-\eta{}^{2^{i-1}}})<\infty$
and $\lim_{k\rightarrow\infty}\|(D^{(k)})^{-\frac{1}{2}}D(D^{(k)})^{-\frac{1}{2}}-I\|_{F}=0$,
\begin{align*}
\|(\overline{D})^{-\frac{1}{2}}D(\overline{D})^{-\frac{1}{2}}-I\|_{F}\leq\sum_{k=1}^{\infty}(\frac{\eta}{1-\eta})^{k} & =\frac{\eta}{1-\eta}\frac{1-\eta}{1-2\eta}=\frac{\eta}{1-2\eta}.
\end{align*}

	Now by concavity of $P(D,\kappa)$ in $D$, we have 
	\begin{align*}
	P(D,\kappa)-P(\overline{D},\kappa) & \leq Tr(\nabla_{D}P(\overline{D},\kappa))\cdot(D-\overline{D})\\
	& =Tr((-\overline{R}^{-1}+\kappa\overline{S}^{-1}+\overline{D}^{-1})\cdot(D-\overline{D}))\\
	& =Tr((-\overline{R}^{-1}+\overline{X}+\kappa\overline{S}^{-1}-\kappa\overline{Y}+\overline{D}^{-1}-\overline{Z})\cdot(D-\overline{D}))\\
	& =Tr((-\overline{R}^{-1}+\overline{X})(\overline{R}-R))+Tr((\overline{S}^{-1}-\overline{Y})(S-\overline{S}))+Tr((\overline{D}^{-1}-\overline{Z})\cdot(D-\overline{D}))
	\end{align*}
	The three trace quantities are all $O(\eta^{2})$. For example,
	\begin{align*}
	Tr((\overline{D}^{-1}-\overline{Z})\cdot(D-\overline{D})) & =Tr((\overline{D})^{\frac{1}{2}}(\overline{D}^{-1}-\overline{Z})(\overline{D})^{\frac{1}{2}}\cdot(\overline{D})^{-\frac{1}{2}}(D-\overline{D})(\overline{D})^{-\frac{1}{2}})\\
	& \leq\|(I-(\overline{D})^{\frac{1}{2}}\overline{Z}(\overline{D})^{\frac{1}{2}})\|_{F}\cdot\|(I-(\overline{D})^{-\frac{1}{2}}D(\overline{D})^{-\frac{1}{2}}\|_{F}\\
	& \leq\eta\cdot\frac{\eta}{1-2\eta}
	\end{align*}
	and similarly for the other two quantities. 
	
	We have thus proved that an $O(\eta)$ approximate analytic center is $O(\eta)$ close to the  exact analytic center, \emph{independent} of $\kappa$. Moreover, the potential functions satisfy $	P(\kappa)-P(\overline{D},\kappa) = O(\eta^2)$.
	
	Now let $\overline{D}^{0}$ be an $O(\delta)$ approximate center for $\Gamma(\kappa^{0})$,
	and $\overline{D}^{1}$ be the approximate center for $\Gamma(\kappa^{1})$
	obtained from $\overline D^0$ using the updates in \eqref{eq:initial-set} followed by a Newton update. Recall that $\kappa^{1}=\kappa^{0}-\frac{\beta}{Tr((\overline{D})^{0}\cdot(\kappa^{0}(\overline{D})^{0}-M)^{-1})}.$
	\begin{align*}
	P(\kappa^{1},\overline{D}^{1})-P(\kappa^{0},\overline{D}^{0}) & =P(\kappa^{1},\overline{D}^{1})-P(\kappa^{1},D^{1})+P(\kappa^{1},D^{1})-P(\kappa^{0},D^{0})+P(\kappa^{0},D^{0})-P(\kappa^{0},\overline{D}^{0})\\
	& \leq-\beta+P(\kappa^{0},D^{0})-P(\kappa^{0},\overline{D}^{0})+P(\kappa^{1},D^{1})-P(\kappa^{1},\overline{D}^{1})\\
	& \leq-\beta+c\cdot\delta^{2}
	\end{align*}
	With $\delta$ and $\beta$ sufficiently small, $-\beta+c\cdot\eta^{2}\leq-c'\beta$
	for some $c'\in(0,1)$.
	\hfill \halmos
\endproof
 \section{First Order Methods for Optimal Diagonal Preconditioning}
\label{sec:first-order-methods}
In this section, we briefly discuss two first order algorithms for finding (right) optimal diagonal preconditioners that do not rely on solving SDPs with Newton updates. They can be easily generalized to the two-sided problem. Their feasibilities have been confirmed in numerical experiments on small matrices, although we leave their customization and optimization on larger problems to future works.

\subsection{A Projected Sub-gradient Descent Algorithm}

To find an optimal (right) diagonal preconditioner for a positive definite
matrix $M=A^TA$, we propose to use projected subgradient descent to solve
the following problem:
\begin{align*}
\min\log(\kappa(DMD))\\
D\succeq I\\
C\succeq D
\end{align*}
 where, for some pre-specified large $C$, the projection step simply
truncates the updated $D_{k}$ at the $k$-th step by an upper bound
$C$ and lower bound $1$. The computation of the subgradient is simple
for this problem, since for the largest eigenvalue $\lambda_{\max}(X)$
of a matrix $X$ and any eigenvector $v$ associated with $\lambda_{\max}(X)$,
$vv^{T}\in\partial\lambda_{\max}(X)$, and similarly $uu^{T}\in\partial(\lambda_{\min}(X))$
for any eigenvector $u$ associated with $\lambda_{\min}(X)$. Thus,
\begin{align*}
(\frac{v(DMD)v^{T}(DMD)}{\lambda_{\max}(DMD)}\times MD)\cdot\mathbf{1}+\mathbf{1}^{T}\cdot(DM\times\frac{v(DMD)v^{T}(DMD)}{\lambda_{\max}(DMD)}) & \in\partial_{D}\log\lambda_{\max}(DMD)\\
(\frac{u(DMD)u^{T}(DMD)}{\lambda_{\min}(DMD)}\times MD)\cdot\mathbf{1}+\mathbf{1}^{T}\cdot(DM\times\frac{u(DMD)u^{T}(DMD)}{\lambda_{\min}(DMD)}) & \in\partial_{D}(\log\lambda_{\min}(DMD))
\end{align*}
where $\times$ is the entry-wise matrix multiplication and $\cdot$
is the usual matrix-vector multiplication. The update at the $k$-th
step is then
\begin{align*}
D_{k}-\alpha_{k}(\mathbf{1}^{T}\cdot D_{k}M\times P_{k}+P_{k}\times MD_{k}\cdot\mathbf{1})\\
P_{k}=\frac{v(DMD)v^{T}(DMD)}{\lambda_{\max}(DMD)}-\frac{u(DMD)u^{T}(DMD)}{\lambda_{\min}(DMD)}
\end{align*}
Extreme eigenvectors and eigenvalues in the subgradient can be efficiently computed with the Lanczos algorithm \citep{lanczos1950iteration}.
For step sizes, we can choose $\alpha_{k}=1/k$ or $\alpha_{k}=1/\sqrt{k}$. The corresponding projected sub-gradient algorithms for left and two-sided preconditionings can be derived similarly.

\subsection{A First-Order Potential Reduction Algorithm}
Recall the potential reduction algorithm introduced in Section \ref{sec:optimal}, which aims to minimize the following potential function:
\begin{align*}
P(\kappa) =\max_{D\succeq0} & \log\det(M-D)+\log\det(\kappa D-M)+\log\det D
\end{align*}
As an alternative to Algorithms \ref{alg:exact} and \ref{alg:approximate} which are based on Newton updates with Nesterov-Todd directions, we suggest to achieve potential reduction using projected alternating gradient descent-ascent on the following min-max objective: 
\begin{align*}
\min_{\kappa\geq1} \max_{D\succeq0} & \log\det(M-D)+\log\det(\kappa D-M)+\log\det D
\end{align*}
using projected updates of $D$ whenever necessary to ensure that the inner maximization problem has finite value. A detailed study of this proposal is left for future work.

 \end{APPENDICES}

\end{document}